\documentclass[11pt]{article}
\usepackage{amsfonts}
\usepackage{amsfonts}
\usepackage{amssymb,latexsym}

\usepackage{amsmath}
\usepackage{amsthm}
\newtheorem{theorem}{Theorem}

\newtheorem{proposition}[theorem]{Proposition}
\newtheorem{lemma}[theorem]{Lemma}

\newtheorem{corollary}[theorem]{Corollary}

\numberwithin{equation}{section}

\newtheorem{Thm}{Theorem}[section]

\begin{document}
\setlength{\baselineskip}{1.2\baselineskip}
\title
{ Convexity of Level Sets of  Minimal Graph  on Space Form with Nonnegative Curvature }
\author{{Peihe Wang$^1$ \thanks {Email: $^1${peihewang@hotmail.com} }} \\
(\small{Sch. of Math. Sci., Qufu Normal University,
273165, Qufu Shandong, China })\\
     Dekai Zhang$^2$\thanks{ E-mail:  $^2$ dekzhang@mail.ustc.edu.cn}\\
     (\small{Sch. of Math. Sci., Univ. of Sci. and Tech. of China,
230026, Hefei Anhui, China })}
\date{}
 \maketitle
{\bf Abstract:}\ For the minimal graph defined on a convex ring in the space form with nonnegative curvature, we obtain the regularity and the strict convexity about its level sets by the continuity method.



\section{\bf Introduction}
\noindent

The Geometry especially the convexity of level sets of the solutions to elliptic
partial differential equations has been interesting to us for a long time. For
instance, Alfhors(\cite{AH}) concluded that
level curves of Green function on simply connected convex domain in
the plane are the convex Jordan curves. Shiffman(\cite{Sh56})
studied the minimal annulus in $R^3$ whose boundary
consists of two closed convex curves in parallel planes $P_1,P_2$,
he derived that the intersection of the surface with any parallel
plane $P$, between $P_1$ and $P_2$, is a convex Jordan curve. In
1957, Gabriel(\cite{Ga57}) proved that the level sets of the Green
function on a 3-dimensional bounded convex domain are strictly
convex and Lewis(\cite{Lew77}) extended Gabriel's result to $p$-harmonic
functions in higher dimensions. Makar-Limanov(\cite{ML71}) and
Brascamp-Lieb(\cite{BL}) got the results on the Poisson
equation and first eigenvalue equation with Dirichlet boundary value
problem on bounded convex domain. Caffarelli-Spruck(\cite{CS82})
generalized Lewis's results(\cite{Lew77}) to a class of semilinear
elliptic partial differential equations. Motivated by the result of
Caffarelli-Friedman(\cite{CF85}), Korevaar(\cite{Korevaar}) gave a new
proof on the results of Gabriel and Lewis (\cite{Ga57},
\cite{Lew77}) using the deformation process and the constant rank
theorem of the second fundamental form of convex level sets of
$p$-harmonic function. Moreover, he also concluded in his paper(\cite{Korevaar}) that level sets of minimal graph defined on convex rings are strictly convex. Kawohl(\cite{Ka85}) gave a  survey of this subject. For more recent related extensions, please see the
papers by Bianchini-Longinetti-Salani(\cite{BLS}), Xu(\cite{XuLu}) and Bian-Guan-Ma-Xu(\cite{BGMX}).

 On the curvature estimates of the level sets, Ortel-Schneider(\cite{OS}),  Longinetti(\cite{Lo83}, \cite{Lo87})proved that the curvature of level curves
  attains its minimum on the boundary (see also Talenti\cite{T83} for related
  results) for  2-dimensional
 harmonic function with convex level curves.
   Furthermore, Longinetti  studied the precise relation between
  the curvature of the convex level lines and the height of minimal graph in \cite{Lo87}. The curvature estimate of the level sets of the solution to partial differential equations then have no new progress until recently, Ma-Ou-Zhang(\cite{MOZ09}) got the
Gaussian curvature estimates of the convex level sets
 of harmonic functions which depend on the Gaussian
curvature of the boundary  and the norm of the gradient on the
boundary in $\mathrm{R}^{n}$. Furthermore, in \cite{MXZ} the concavity of the Gaussian curvature of the convex level sets of
$p$-harmonic functions with respect to the height was derived to describe the variation of the curvature along the height of the function. In \cite{GX}, the lower bound of the principal curvature of the convex level sets of the solution to a kind of fully nonlinear elliptic equations was derived. For Poisson equations and a class of semilinear
elliptic partial differential equations, Caffarelli-Spruck(\cite{CS82}) concluded that the level sets of their solutions are all convex with respect to the gradient direction, the curvature estimate of the level sets has been got by Wang-Zhang(\cite{WZ}), and in the same paper they also described the geometrical properties of the level sets of the minimal graph. In the sequel, following the technique in \cite{MXZ}, Wang(\cite {WangPJM}) got the precise relation between
  the curvature of the convex level  sets and the height of minimal graph of general dimensions which generalized the previous results of Longinetti(\cite{Lo87}).

For the Riemannian manifold case, Papadimitrakis(\cite{papad}) concluded the convexity of the level curves of harmonic functions on convex rings in the hyperbolic plane via one complex variable tools. Ma-Zhang(\cite{Mzyb}) generalized Papadimitrakis's results to space form of general dimensions. Partial results in \cite{Mzyb} can be stated as follows.
\begin{theorem}\label{MaZhang}({\bf\cite{Mzyb}})
Let $(M^n,g)$ be a space form with constant sectional curvature $1$ or $-1$, and $\Omega_0$ and $\Omega_1$ be bounded smooth strictly convex domains in $M^n$, $n\ge 2$ and $\bar{\Omega}_1\subseteq \Omega_0$. Let $\omega$ satisfy
\begin{align*}
\left\{ {\begin{array}{lll}
\triangle \omega=0 & \mathrm{\mbox {in}} & \ \Omega=\Omega_0\setminus\bar{\Omega}_1,\\
\omega=0& \mathrm{\mbox {on}}& \  \partial \Omega_0,\\
\omega=1&\mathrm{ \mbox {on}}&  \ \partial \Omega_1.
\end{array}} \right.
\end{align*}
Then $\nabla \omega\neq 0$ is valid everywhere in $\Omega$ and all the level sets of $\omega$ are strictly convex with respect to $\nabla \omega$.
\end{theorem}

Based on the above strict convexity of the level sets of harmonic funcions defined on the convex ring in space forms, following the technique in \cite{Korevaar}, we come in this paper to consider another important geometrical object, the minimal graph defined on the convex ring in space form with nonnegative curvature. We mainly get the following theorem.

\begin{theorem}\label{WZ}
Let $(M^n,g)$ be a space form with constant sectional curvature $\epsilon\ge0$ , and $\Omega_0$ and $\Omega_1$ be bounded smooth strictly convex domains in $M^n$, $n\ge 2$ and $\bar{\Omega}_1\subseteq \Omega_0$. Then the following minimal graph equation defined on $\Omega=\Omega_0\setminus\bar{\Omega}_1$
\begin{align}\label{minimalgraph}
\begin{split}
\left\{ {\begin{array}{lll}
\mathrm{div}(\frac{\nabla u}{\sqrt{1+|\nabla u|^2}}) =0 & \mathrm{\mbox {in}} & \ \Omega=\Omega_0\setminus\bar{\Omega}_1,\\
u=0& \mathrm{\mbox {on}}& \  \partial \Omega_0,\\
u=1& \mathrm{\mbox {on}}&  \ \partial \Omega_1
\end{array}} \right.
\end{split}
\end{align}
has a unique smooth solution $u$. Moreover, $\nabla u \neq 0$ is valid everywhere in $\Omega$ and all the level sets of $u$ are strictly convex with respect to $\nabla u$.
\end{theorem}

Following this theorem, we immediately get a geometrical property of the minimal graph defined on the convex ring.
\begin{corollary}\label{WZCor}
Let $(M^n,g)$ be a space form with constant sectional curvature $\epsilon\ge0$ , and $\Omega_0$ and $\Omega_1$ be bounded smooth strictly convex domains in $M^n$, $n\ge 2$ and $\bar{\Omega}_1\subseteq \Omega_0$. Let $u$ be the solution to \eqref{minimalgraph}. Then $|\nabla u|$ increases strictly along the gradient direction.
\end{corollary}

Remark that the case of $\epsilon =0$ has been already concluded in \cite{Korevaar}. And it is a pity that we can not derive similar information or counterexample to the case $\epsilon<0.$  Note that in \cite{YingShih} a counterexample was constructed to show that the level sets of the first eigenfunction of a domain with negative curvature are not convex.

The paper is organized as follows: in section 2, we list the notations and the preliminaries being used during the process of the proof. In section 3, we give the existence of the solution of minimal graph defined on convex rings in space form and deduce the $C^{2,\alpha}$ continuity of the solution with respect to the boundary value. In section 4, we prove a constant rank theorem of the second fundamental form of the level sets of minimal graph. In section 5, we prove the regularity and strict convexity of the level sets of minimal graph.

\section{Notations and Preliminaries}
\setcounter{equation}{0} \setcounter{theorem}{0}
\noindent

 In this section, we introduce some notations and preliminaries for the main results.

Firstly, we give the derivative commutation formula in Riemannian geometry.

\begin{lemma}\label{Derivative comutation}
Let $u$ be a smooth function defined on Riemannian manifold $M^n$ with constant sectional curvature and denote by $R_{ijkl}$ the curvature tensor of $M^n$, then
\begin{align*}
u_{ijk}=u_{ikj}+u_lR_{lijk}
\end{align*}
and
\begin{align*}
u_{ijkl}=u_{klij}-u_{\xi j}R_{\xi kl i}-u_{ i\xi}R_{\xi kl j}+u_{\xi l}R_{\xi ijk}+u_{k\xi}R_{\xi ijl}.
\end{align*}
Remark that Einstein summation convention here is adopted.\hfill $\sharp$
\end{lemma}

In the sequel, we list out a linear algebra formula(\cite{ChHu, CMSalani, HuM}) used frequently in our proof.

\begin{proposition}\label{SigmaK}
Let $u$ be a smooth function and $(D^2u)$ be its Hessian matrix. Assume that
$$(D^2u)=\left( {\begin{array}{*{20}c}
   {u_{ij} } & {u_{in} }  \\
   {u_{nj} } & {u_{nn} }  \\
\end{array}} \right)_{n \times n}
.$$
Denote by $\sigma_k(A)$ the $k-$th elementary symmetric function of the eigenvalue of the matrix $A$. Then we have
\begin{align}\label {Formula}
\begin{split}
\sigma_{l+1}(D^2u)=&\sigma_{l+1}(u_{ij})+u_{nn}\sigma_{l}(u_{ij})-\sum\limits_{i}u_{ni}u_{in}\sigma_{l-1}(u_{pq}|i)\\
&+\sum\limits_{i\neq j}u_{ni}u_{jn}u_{ij}\sigma_{l-2}(u_{pq}|ij)-\sum\limits_{i\neq j,j\neq k,k\neq i}u_{ni}u_{jn}u_{ik}u_{kj}\sigma_{l-3}(u_{pq}|ijk)+T,
\end{split}
\end{align}
where we denote by $(u_{pq}|i)$ the symmetric matrix obtained from $(u_{pq})$ by deleting the $i-$row and
$i-$column and by$(u_{pq}|ij)$ the symmetric matrix obtained from $(u_{pq})$ when deleting the $i, j-$rows and $i, j-$columns,
and similarly we define $(u_{pq}|ijk)$, and also every term in the polynomial $T$ includes at least $3$ factors like $u_{rs}$ with $r\neq s$. So if the matrix $(u_{pq})$ is diagonal at a fixed point, we then have at this point that
$$T=0,\  DT=0,\  \mbox{and}\ \  D^2T=0.$$
Remark that all the Latin indices will vary from $1$ to $n-1$.\hfill $\sharp$
\end{proposition}

For the calculation principle of $\sigma_k(A)$, we do not mention much more here and one can refer to \cite{XuLu} for details.

For a smooth function $u$ defined on a Riemannian manifold $M^n$, its graph can be considered as a hypersurface in $M^{n}\times{\mathrm{R}}$ with canonical product Riemannian metric. In \cite{Spruck}, the mean curvature of this hypersurface has been already deduced.

\begin{proposition} \label{Mean curvature of graph}{\bf(\cite{Spruck})}
Let $u:M^{n}\rightarrow {\mathrm{R}}$ be a smooth function defined on a Riemannian manifold $M^n$. Consider the graph of $u$, denoted by $\Sigma_u=F(M^{n})$,  where $F:M^{n}\rightarrow{M^{n}\times{\mathrm{R}}}$ is defined to be $F(p)=(p,\,u(p))$ and $M^{n}\times{\mathrm{R}}$ is equipped with the canonical product Riemannian metric. Then the mean curvature of $\Sigma_u$ is
$$H=-\mathrm{div}(\frac{\nabla u}{\sqrt{1+|\nabla u|^2}}),$$
here, $\mathrm{div}$ is the divergence operator on $M^n$.\hfill $\sharp$

\end{proposition}

The following corollary is obvious.
\begin{corollary} \label{minimal graph equation}
Let $u$ be a smooth function defined on  $M^{n}$ , then the graph of $u$, $\Sigma_u$, is minimal in $M^{n}\times{\mathrm{R}}$ if and only if $\mathrm{div}(\frac{\nabla u}{\sqrt{1+|\nabla u|^2}})=0$ holds on $M^{n}$.\hfill $\sharp$
\end{corollary}

We remark that in \cite{Colding} the equation of minimal graph defined on manifolds is also deduced.

The following lemma is essential to express the geometric quantity of the level sets in terms of the function itself. We omit the proof and one can refer to Proposition 2.1 in \cite{Mzyb}.
\begin{lemma} \label{AAA}
Let $u$ be a smooth function defined on  $M^{n}$ with nonzero gradient everywhere. Assume that the level sets of $u$ is convex with respect to the normal direction. Let $\{e_\alpha,\ \alpha=1,2,\cdots,\ n\}$ be a local orthogonal frame on $M^{n}$. Then the $k-$th curvature of the level set $\Sigma^c=u^{-1}(c)$ is
\begin{align} \label{ k-th curvature}
\begin{split}
\sigma_k[\Sigma^c]=(-1)^k\sum_{\alpha,\beta=1}^n\frac{\partial \sigma_{k+1}(D^2u)}{\partial u_{\alpha\beta}}u_\alpha u_\beta|\nabla u|^{-(k+2)},
\end{split}
\end{align}
where $1\le k\le n-1$.

 Furthermore, if we take $e_n=\frac{\nabla u}{|\nabla u|}$ as the unit normal direction, then the second fundamental form of the level set of $u$ is
\begin{align} \label{BBB}
\begin{split}
h_{ij}=-\frac{u_{ij}}{|\nabla u|}.
\end{split}
\end{align}

\end{lemma}

\section{\bf  Existence of the Solution}
\setcounter{equation}{0} \setcounter{theorem}{0}
\noindent

In this section, we settle the existence of the solution to minimal graph equation and deduce some properties of it.

Firstly, we construct a supersolution to the minimal graph equation.

Given a smooth and strictly convex  ring $\Omega=\Omega_0\backslash\Omega_1$ in a space form. According to Theorem \ref{MaZhang}, there exists a unique harmonic function $\omega$ defined on $\Omega$ such that $\omega=0,\tau(0<\tau\le 1)$ on $\partial \Omega_0,\ \partial \Omega_1$ respectively. We will construct a supersolution in terms of this harmonic function.

Let $g(t)=-\frac{t^2}{4\tau}+\frac{5t}{4}$ and $v=g(\omega)$, we conclude that $v$ is a supersolution of the following minimal graph equation
 \begin{align}\label {boudary-tau}
\begin{split}
\left\{ {\begin{array}{lll}
\mathrm{div}(\frac{\nabla u}{\sqrt{1+|\nabla u|^2}}) =0 & \mathrm{\mbox {in}} & \ \Omega=\Omega_0\setminus\bar{\Omega}_1,\\
u=0& \mathrm{\mbox {on}}& \  \partial \Omega_0,\\
u=\tau& \mathrm{\mbox {on}}&  \ \partial \Omega_1\,.
\end{array}} \right.
\end{split}
\end{align}
 In fact, we note by the strict convexity of the level sets of $\omega$ and Lemma \ref{AAA} that $|\nabla \omega|>0$ and $\omega_{ii}<0$ for $i=1,\,2,\,\cdots,n-1$. Thus, for the suitable frame such that $\omega_n=\frac{\nabla \omega}{|\nabla \omega|}$ and $\omega_{i}=0$ for $i=1,\,2,\,\cdots,n-1$, we have
\begin{align}\label {}
\begin{split}
\frac{{\sum\limits_{\alpha ,\beta  = 1}^n {\omega _{\alpha \beta } \omega _\alpha  \omega _\beta  } }}{{|\nabla \omega |^2 }} = \omega _{nn}  = \Delta \omega  - \sum\limits_{i = 1}^{n - 1} {\omega _{ii} }  = 0 - \sum\limits_{i = 1}^{n - 1} {\omega _{ii} }  > 0.
\end{split}
\end{align}
It is a direct observation that
$$
g(0) = 0,\ g(\tau) = \tau,\ \frac{3}{4}\le g'(t) \le \frac{5}{4} (t \in [0,\tau]),\ g''(t) =  - \frac{1}{2\tau}.$$
Obviously, $v=0 \ \mbox{on}\  \partial \Omega_0$ and $v=\tau \ \mbox{on} \ \partial \Omega_1.$
It also follows that
\begin{align}\label {}
\begin{split}
v_\alpha   =& g'\omega _\alpha ,\ |\nabla v|^2  = (g')^2 |\nabla \omega |^2 ,\ v_\alpha  v_\beta   = (g')^2 \omega _\alpha  \omega _\beta  , \\
 v_{\alpha \beta }  =& g''\omega _\alpha  \omega _\beta   + g'\omega _{\alpha \beta }  .
\end{split}
\end{align}
Therefore, using the Einstein summation convention,
\begin{align}\label {Lv}
\begin{split}
 {\rm{L}}v = &[(1 + |\nabla v|^2 )\delta _{\alpha \beta }  - v_\alpha  v_\beta  ]v_{\alpha \beta }  \\
  =& [1 + (g')^2 |\nabla \omega |^2 ](g''|\nabla \omega |^2  + g'\Delta \omega ) - (g')^2 \omega _\alpha  \omega _\beta  (g''\omega _\alpha  \omega _\beta   + g'\omega _{\alpha \beta } ) \\
  =& [1 + (g')^2 |\nabla \omega |^2 ](g''|\nabla \omega |^2 ) - (g')^2 \omega _\alpha  \omega _\beta  (g''\omega _\alpha  \omega _\beta   + g'\omega _{\alpha \beta } ) \\
  =& g''|\nabla \omega |^2  - (g')^3 \omega _\alpha  \omega _\beta  \omega _{\alpha \beta }  \\
  \le&  - \frac{1}{2\tau}|\nabla \omega |^2  < 0.
\end{split}
\end{align}
Thus we have constructed a supersolution to the minimal graph equation defined on a convex ring in space forms.

The following proposition will show that the maximum of the norm of gradient of minimal graph could be attained on the boundary. More precisely,
\begin{proposition} \label{Norm of Gradient}
Let $u$: $M\rightarrow \mathrm{R}$ be a minimal graph defined on $\Omega$ in an $n-$dimensional  space form $M$ with nonnegative curvature. Then we have
\begin{align}\label {NormEstimate}
\begin{split}
\sup_{x\in\Omega}|\nabla u|(x) \le \sup_{x\in\partial\Omega}|\nabla u|(x).
\end{split}
\end{align}
\end{proposition}

{\bf Proof:\ \ } Let $\phi=\frac{1}{2}|\nabla u|^2$. We choose the frame $\{e_\alpha\}$ on the manifold such that the Riemannian curvature tensor takes the form
$R_{\alpha\xi\beta\eta}=\epsilon(\delta_{\alpha\beta}\delta_{\xi\eta}-\delta_{\alpha\eta}\delta_{\beta\xi})$, where $\epsilon\ge 0$ and $\delta_{\alpha\beta}$ is the Kronecker symbol.
It follows that
\begin{align}\label {Gradient}
\begin{split}
\sum\limits_{\alpha ,\beta }& {[(1 + |\nabla u|^2 )\delta _{\alpha \beta }  - u_\alpha  u_\beta  ]\phi _{\alpha \beta } }\\
  =& \sum\limits_{\alpha ,\beta ,\gamma } {[(1 + |\nabla u|^2 )\delta _{\alpha \beta }  - u_\alpha  u_\beta  ](u_\gamma  u_{\gamma \alpha \beta } {\rm{ + }}u_{\gamma \alpha } u_{\gamma \beta } )}  \\
  =& \sum\limits_{\alpha ,\beta ,\gamma } {u_\gamma  [(1 + |\nabla u|^2 )\delta _{\alpha \beta }  - u_\alpha  u_\beta  ]u_{\gamma \alpha \beta } } {\rm{ + }}(1 + |\nabla u|^2 )\sum\limits_{\alpha ,\beta } {u_{\alpha \beta } ^2 }  - \sum\limits_\gamma  {\phi _\gamma  } ^2 \ ,
\end{split}
\end{align}
by Lemma \ref{Derivative comutation} we have
\begin{align}\label {}
\begin{split}
\sum\limits_{\alpha ,\beta } &{[(1 + |\nabla u|^2 )\delta _{\alpha \beta }  - u_\alpha  u_\beta  ]\phi _{\alpha \beta } }  \\
 \ge & \sum\limits_{\alpha ,\beta ,\gamma } {u_\gamma  [(1 + |\nabla u|^2 )\delta _{\alpha \beta }  - u_\alpha  u_\beta  ](u_{\alpha \beta \gamma }  + u_\xi  R_{\zeta \alpha \gamma \beta } )}  - \sum\limits_\gamma  {\phi _\gamma  } ^2\\
  \ge &  - \sum\limits_{\alpha ,\beta ,\gamma } {u_\gamma  [(1 + |\nabla u|^2 )\delta _{\alpha \beta }  - u_\alpha  u_\beta  ]_\gamma  u_{\alpha \beta } }  - \sum\limits_\gamma  {\phi _\gamma  } ^2  \\
  = & - 2\Delta u\sum\limits_\gamma  {u_\gamma  \phi _\gamma  }  + 2\sum\limits_\gamma  {\phi _\gamma  } ^2  - \sum\limits_\gamma  {\phi _\gamma  } ^2  \\
  = &\sum\limits_\gamma  {\phi _\gamma  } ^2  - 2\Delta u\sum\limits_\gamma  {u_\gamma  \phi _\gamma  }  \ .
\end{split}
\end{align}
Thus by the weak maximum principle we reach the conclusion.\hfill $\sharp$

Now, we can get the existence of the minimal graph on the convex ring in space forms with nonnegative curvature.
\begin{theorem}\label{WZtau}
Let $(M^n,g)$ be a space form with constant sectional curvature $\epsilon\ge0$ , and $\Omega_0$ and $\Omega_1$ be bounded smooth strictly convex domains in $M^n$, $n\ge 2$ and $\bar{\Omega}_1\subseteq \Omega_0$. Then the following minimal graph equation defined on $\Omega$ for $0<\tau \le 1$
\begin{align*}
\left\{ {\begin{array}{lll}
\mathrm{div}(\frac{\nabla u^\tau}{\sqrt{1+|\nabla u^\tau|^2}}) =0 & \mathrm{\mbox {in}} & \ \Omega=\Omega_0\setminus\bar{\Omega}_1,\\
u^\tau=0& \mathrm{\mbox {on}}& \  \partial \Omega_0,\\
u^\tau=\tau& \mathrm{\mbox {on}}&  \ \partial \Omega_1
\end{array}} \right.
\end{align*}
has a unique smooth solution $u^\tau$. Moreover,

$\mathrm{(i)}$. there exists a positive constant $C_1>0$ such that $\sup\limits_{x\in \bar{\Omega}}|\nabla u^\tau|\le C_1\tau\,,$

$\mathrm{(ii)}$. there exists a positive constant $C_2>0$ such that $$||u^t- u^\tau||_{C^{2,\alpha}(\bar{\Omega})} \le C_2|t-\tau|$$ for some $\alpha \in (0,\,1)$ and any $t\in (0,\,1]$.
 Especially, setting $t\rightarrow 0$, we get that $$|| u^\tau||_{C^{2,\alpha}(\bar{\Omega})} \le C_2\tau.$$
\end{theorem}

{\bf Proof:\ \ } According to Propositon \ref{Norm of Gradient} and the supersolution we constructed, we then deduce that the solution to the minimal graph equation(\ref{boudary-tau}) has a priori gradient estimate as follows
\begin{align}\label {GradientUpperBound}
\begin{split}
\sup_{x\in\bar{\Omega}}|\nabla u|(x) & \le \sup_{x\in\partial\Omega}|\nabla u|(x) \le \sup_{x\in\partial\Omega}|\nabla v|(x)\\
&\le\sup_{x\in \partial\Omega_0 \bigcup  \partial\Omega_1} \left( g'(\omega)|\nabla \omega|\right)(x)\le \frac{5}{4}\sup_{x\in \partial\Omega_1}|\nabla \omega|(x).
\end{split}
\end{align}
Note that the maximum of $|\nabla \omega|$ can only be attained on the interior boundary $\partial\Omega_1$, one can refer to Proposition 4.1 in \cite{MOZ09}.

Now, by the classical theory of quasilinear elliptic partial differential eqations(\cite{Tru}: Ch.13 and Ch.16) we then get the existence theorem.

Uniqueness and smoothness of the solution are obvious.

For $\mathrm{(i)}$, by \eqref{GradientUpperBound}, we only need to consider the gradient estimate of $\omega^\tau$, the harmonic function defined on the same convex rings and with the same boundary value, i.e.
\begin{align}\label {omegatau}
\begin{split}
\left\{ {\begin{array}{lll}
\triangle  \omega^\tau =0 & \mathrm{\mbox {in}} & \ \Omega=\Omega_0\setminus\bar{\Omega}_1,\\
\omega^\tau=0& \mathrm{\mbox {on}}& \  \partial \Omega_0,\\
\omega^\tau=\tau& \mathrm{\mbox {on}}&  \ \partial \Omega_1\,.
\end{array}} \right.
\end{split}
\end{align}

It is obvious that $\omega^\tau=\tau\omega^1$ and it follows that $|\nabla \omega^\tau |(x)= \tau |\nabla \omega^1 |(x)$ for $x\in \Omega$, this gives the proof of (1). Additionally, since $|\nabla \omega^1|\neq 0 $ holds on $\bar{\Omega}$, we then can take a positive constant $C_0>0$ such that for $x\in \Omega$,
\begin{align}\label {lip0}
\begin{split}
| \nabla \omega^\tau |(x) \geq C_0 \tau\,.
\end{split}
\end{align}

For $\mathrm{(ii)}$, we need more than a word. Firstly, we note that
 \begin{align}\label {lip2}
\begin{split}
|| \omega^\tau ||_{C^{2,\alpha}(\bar{\Omega})}= \tau || \omega^1 ||_{C^{2,\alpha}(\bar{\Omega})}\leq C \tau
\end{split}
\end{align}
for some $\alpha \in (0,\,1)$ and $C>0$.

Without loss of generality we assume that $t>\tau$. Let $h=u^t-u^\tau$, then $h$ would satisfy a linear divergence type equation. In fact, $h$ has boundary value $0$, $t-\tau$ on the boundary $\partial \Omega_0$, $\partial \Omega_1$, respectively. Furthermore, if we rewrite the minimal graph equation to the form $\mathrm{div}(A(\nabla u))=0$, we have by following the discussion in \cite{John} or the theory in \cite{Tru} that
 \begin{align}\label {Lipcchtiz-1}
\begin{split}
0=&\mathrm{div}(A(\nabla u^t))-\mathrm{div}(A(\nabla u^\tau))=\sum\limits_\alpha D_\alpha\left(A^\alpha(\nabla u^t)-A^\alpha(\nabla u^\tau)\right)\\
=&\sum\limits_\alpha D_\alpha\left(A^\alpha(s\nabla u^t+(1-s)\nabla u^\tau)|^1_0\right)\\
=&\sum\limits_\alpha D_\alpha\left(m_{\alpha\beta}(x)D_\beta h\right)\ ,
\end{split}
\end{align}
where,
$$ m_{\alpha\beta}(x) =\int_0^1\frac{\partial A^\alpha}{\partial p_\beta}\left(s\nabla u^t+(1-s)\nabla u^\tau\right)\mathrm{d}s\,.$$

Anyway, $h(x)$ satisfies the following uniformly elliptic differential equation for the uniform gradient bound we have deduced,
\begin{align*}
\left\{ {\begin{array}{lll}
D_\alpha\left(m_{\alpha\beta}(x)D_\beta h\right) =0 & \mathrm{\mbox {in}} & \ \Omega=\Omega_0\setminus\bar{\Omega}_1,\\
h=0& \mathrm{\mbox {on}}& \  \partial \Omega_0,\\
h=t-\tau& \mathrm{\mbox {on}}&  \ \partial \Omega_1\,.
\end{array}} \right.
\end{align*}

In order to give the $C^{2,\alpha}$ estimate of $h$, according to Theorem 8.33 of \cite{Tru} or Theorem $\mathrm{8.33}'$ in \cite{John}, we need to extend the boundary value of $h$ to a smooth function defined on the whole domain $\bar{\Omega}$. Fortunately, the harmonic function $\omega^{t-\tau}$ is just suitable for this target. Thus we have by the maximum principle and \eqref{lip2}
 \begin{align}\label {}
\begin{split}
||h||_{C^{2,\alpha}(\bar{\Omega})}\le C \left( ||h||_{L^\infty(\Omega)}+||\omega^{t-\tau}||_{C^{2,\alpha}(\bar{\Omega})}    \right)\le C_2(t-\tau).
\end{split}
\end{align}
  This completes the whole proof of the theorem.\hfill $\sharp$

\section{\bf A Constant Rank Theorem }
\setcounter{equation}{0} \setcounter{theorem}{0}
\noindent
The constant rank theorem joined with deformation process is usually applied to prove the convexity of the solution or the convexity of the level sets of the solution. Based on this essential tool, lots of important results concerning to convexity appeared(\cite{CGM2007,GLM,GM2003,GMZ,Mzyb}) recently. In this section, we will show a constant rank theorem as follows.
\begin{Thm}\label{WZ3}
Let $\Omega$ be a smooth bounded connected domain on the space form $M^{n}$ with constant curvature $\epsilon \ge 0$. Let $u\in{C^{4}(\Omega)\bigcap{C^{2}(\overline{\Omega})}}$ be the solution to the prescribed mean curvature equation
$$\mathrm{div}(\frac{\nabla u}{\sqrt{1+|\nabla u|^2}})=H(x),$$
where $H(x)\ge0$ satisfies the structure condition
\begin{align}\label{structurecondition}
\begin{split}
3H_\alpha H_\beta+4\epsilon H^2\delta_{\alpha\beta}\le 2HH_{\alpha\beta}.
\end{split}
\end{align}
 Assume $|\nabla{u}|\neq{0}$ in $\Omega$ and the level sets of $u$ are all convex with respect to the normal $\nabla u$, then the second fundamental form of the level sets of $u$ must have the same rank at all points in $\Omega$.
\end{Thm}

{\bf Proof:\ } We firstly represent the equation into the following
$$\sum\limits_{\alpha,\beta}\left[(1+|\nabla{u}|^{2})\delta_{\alpha\beta}-u_\alpha u_\beta\right]u_{\alpha\beta}=H(x)(1+|\nabla u|^2)^{\frac 3 2}\triangleq f(x,\,\nabla u) \ \  \mathrm{in} \ \ \Omega.$$

The theorem is obviously of local feature, so we may assume that the level set $\Sigma^c=u^{-1}(c)$ is connected for each $c$ in a neighborhood of some $c_0=u(x_0)$. We will compute in a neighborhood of some point $x_0\in \Omega$ such that the minimal rank $l$ of the second fundamental form of $\Sigma^{u(x_0)}$ is attained at $x_0$. Without loss of generality, we can assume that $0\le l \le n-2$, otherwise we then have got the constant rank $n-1$.

Let $U$ be a small open neighborhood of $x_0$ such that for each $x\in U$, there are $l$ ``good" eigenvalues of the second fundamental form of $\Sigma^c$ which are bounded from below by a positive constant, and the other $n-1-l$ eigenvalues of the second fundamental form of $\Sigma^c$ are very small and we name them ``bad" eigenvalues. In the following, we will denote by $G$, $B$ the index set of ``good" eigenvalues and ``bad" eigenvalues respectively.

For any $x\in U$ fixed, we can choose a frame $e_1,\,e_2,\cdots,e_{n-1},\,e_n$ such that
$|\nabla u|(x)=u_n(x)>0$ and the matrix $(u_{ij})(i,\,j=1,2,\cdots,n-1)$ is diagonal at $x$.

By \eqref{BBB}, the second fundamental form of $\Sigma^c$ then is also diagonal at $x$. Without loss of generality we could assume $h_{11}\le h_{22}\le\cdots\le h_{n-1n-1}$ and there exists a positive constant $C>0$ depending only on $||u||_{C^4}$ and the domain $U$ such that $h_{n-ln-l}\ge C$ for all $x\in U$. For convenience, we let $G=\{n-l,\,n-l+1,\,\cdots,\,n-1\}$ and $B=\{1,\,2,\,\cdots,\,n-l-1\}$ be the ``good"  and ``bad" sets of indices, respectively. We also denote
$$B=\{h_{11},\,h_{22},\,\cdots,\,h_{n-l-1n-l-1}\}\ \ \mbox{and}\ \ G=\{h_{n-ln-l},,\,\cdots,\,h_{n-1n-1}\}.$$
Note that for a fixed $\delta >0$, we can choose the neighborhood $U$ small enough such that $h_{jj}<\delta$ for all $j\in B$ and $x\in U$.

Denote by $\lambda_i,\ i=1,2,\cdots,n-1$ the principal curvature of level sets of $u$ and $\lambda=(\lambda_1,\cdots,\lambda_{n-1})$, it is obvious that $\lambda_i=h_{ii}$ under the assumption above. We set $\mu_i=u_{ii},\ \mu=(\mu_1,\cdots,\mu_{n-1})$ during the whole proof, then $\mu_i=-u_n\lambda_i$ for $i=1,2,\cdots,n-1$ at $x$.

Generally, the following notations are usually necessary to conclude the constant rank theorems.

For two functions defined on $\Omega$, $h$ and $k$, for $y\in \Omega$, we denote by $ h(y)\preceq k(y)$ if and only if
$$ (h-k)(y)\le (C_3 \phi +C_4|\nabla \phi|)(y).$$
Also, we call that $h(y)\sim k(y)$ if and only if $h(y)\preceq k(y)$ and $k(y)\preceq h(y)$. Generally, we say that $h\sim k$ is valid if and only if $h(y)\sim k(y)$ for any $y\in\Omega$ holds.

During the whole proof, the Greek indices such as $\alpha,\beta$ etc. will vary from $1$ to $n$ while the Latin indices such as $i,j,k$ etc. will vary from $1$ to $n-1$.

For each $c$, set
\begin{align}\label{varphi}
\begin{split}
\varphi  = |\nabla u|^{l + 3} \sigma _{l + 1} \left( \lambda  \right) = ( - 1)^{l + 1} \sum\limits_{\alpha ,\beta } {\frac{{\partial \sigma _{l + 2} (D^2 u)}}{{\partial u_{\alpha \beta } }}u_\alpha  u_\beta  ,}
\end{split}
\end{align}

From now on, all the calculations will be done at $x$ with the above frame. In order to simplify the process of the proof, we introduce the notation(\cite{Mzyb}) $b_{ij,\xi}$ as follows,
\begin{align}\label{bij-xi}
\begin{split}
u_nu_{ij\xi}=-u_n^2b_{ij,\xi}+u_{ni}u_{j\xi}+u_{nj}u_{i\xi}+u_{n\xi}u_{ij}.
\end{split}
\end{align}

Easy to get for $j\in B$,
$$\lambda_j\sim 0.$$

Now we come to compute the first order derivative of $\varphi$.
\begin{align*}
\varphi _\xi   =& ( - 1)^{l + 1} \left( {\sum\limits_{\alpha ,\beta } {\frac{{\partial \sigma _{l + 2} (D^2 u)}}{{\partial u_{\alpha \beta } }}u_\alpha  u_\beta  } } \right)_\xi   \\
  =& ( - 1)^{l + 1} \sum\limits_{\alpha ,\beta ,\gamma ,\delta } {\frac{{\partial ^2 \sigma _{l + 2} (D^2 u)}}{{\partial u_{\alpha \beta } \partial u_{\gamma \delta } }}u_\alpha  u_\beta  } u_{\gamma \delta \xi }  + ( - 1)^{l + 1} \sum\limits_{\alpha ,\beta } {\frac{{\partial \sigma _{l + 2} (D^2 u)}}{{\partial u_{\alpha \beta } }}\left( {u_{\alpha \xi } u_\beta   + u_\alpha  u_{\beta \xi } } \right)}  \\
   = & O_1  + O_2 ;
\end{align*}

For the term $O_1$, by \eqref{Formula},
\begin{align}\label {O1}
\begin{split}
 O_1  =& ( - 1)^{l + 1} \sum\limits_{\alpha ,\beta ,\gamma ,\delta } {\frac{{\partial ^2 \sigma _{l + 2} (D^2 u)}}{{\partial u_{\alpha \beta } \partial u_{\gamma \delta } }}u_\alpha  u_\beta  } u_{\gamma \delta \xi }
 = ( - 1)^{l + 1} u_n ^2 \sum\limits_{\gamma ,\delta } {\frac{{\partial ^2 \sigma _{l + 2} (D^2 u)}}{{\partial u_{nn} \partial u_{\gamma \delta } }}} u_{\gamma \delta \xi }  \\
  =& ( - 1)^{l + 1} u_n ^2 \sum\limits_i {\frac{{\partial \sigma _{l + 1} (u_{pq} )}}{{\partial u_{ii} }}} u_{ii\xi }
    =( - 1)^{l + 1} u_n ^2 \sum\limits_i {\sigma _l \left( {\mu \left| i \right.} \right)} u_{ii\xi } \\
     \sim & ( - 1)^{l + 1} u_n ^2 \sum\limits_{i \in B} {\sigma _l \left( {\mu \left| i \right.} \right)} u_{ii\xi }
  =  - u_n ^{l + 2} \sum\limits_{i \in B} {\sigma _l \left( {\lambda \left| i \right.} \right)} u_{ii\xi }
  \sim  - u_n ^{l + 2} \sigma _l \left( G \right)\sum\limits_{j \in B} {u_{jj\xi } } ;
 \end{split}
\end{align}

For the term $O_2$,
\begin{align}\label {O2}
\begin{split}
 O_2  =& ( - 1)^{l + 1} \sum\limits_{\alpha ,\beta } {\frac{{\partial \sigma _{l + 2} (D^2 u)}}{{\partial u_{\alpha \beta } }}\left( {u_{\alpha \xi } u_\beta   + u_\alpha  u_{\beta \xi } } \right)}  = 2( - 1)^{l + 1} u_n \sum\limits_\beta  {\frac{{\partial \sigma _{l + 2} (D^2 u)}}{{\partial u_{n\beta } }}u_{\beta\xi } }  \\
  =& 2( - 1)^{l + 1} u_n \frac{{\partial \sigma _{l + 2} (D^2 u)}}{{\partial u_{nn} }}u_{n\xi }  + 2( - 1)^{l + 1} u_n \sum\limits_i {\frac{{\partial \sigma _{l + 2} (D^2 u)}}{{\partial u_{ni} }}u_{i\xi } }  \\
  =& 2( - 1)^{l + 1} u_n \sigma _{l + 1} \left( \mu  \right)u_{n\xi }  + 2( - 1)^l u_n \sum\limits_i {\sigma _l \left( {\mu \left| i \right.} \right)u_{in} u_{i\xi } }  \\
 \sim & 2u_n ^{l + 1} \sum\limits_{i \in B} {\sigma _l \left( {\lambda \left| i \right.} \right)u_{in} u_{i\xi } }
  \sim  2u_n ^{l + 1} \sigma _l \left( G \right)\sum\limits_{j\in B} {u_{jn} u_{j\xi } } ;
\end{split}
\end{align}

Therefore, by \eqref{bij-xi}, \eqref{O1} and \eqref{O2}, we have

\begin{align}\label {varphi-xi}
\begin{split}
 \varphi _\xi   \sim &  - u_n ^{l + 2} \sigma _l \left( G \right)\sum\limits_{j \in B} {u_{jj\xi } }  + 2u_n ^{l + 1} \sigma _l \left( G \right)\sum\limits_{j \in B} {u_{jn} u_{j\xi } }  \\
  = & - u_n ^{l + 1} \sigma _l \left( G \right)\left( {u_n \sum\limits_{j \in B} {u_{jj\xi } }  - 2\sum\limits_{j \in B} {u_{jn} u_{j\xi } } } \right) \\
 \sim & u_n ^{l + 2} \sigma _l \left( G \right)\sum\limits_{j \in B} {b_{jj,\xi } }\,.
\end{split}
\end{align}

So we deduce that
\begin{align}\label {1-ordercondition}
\begin{split}
  \sum\limits_{j \in B} {b_{jj,\xi } }\sim 0,\ \forall\, 1\le \xi\le n.
\end{split}
\end{align}

Now we set out to compute the second order derivative of $\varphi$.
\begin{align*}
\varphi _{\xi \eta }  = &( - 1)^{l + 1} \sum\limits_{\alpha ,\beta } {\left( {\frac{{\partial ^2 \sigma _{l + 2} (D^2 u)}}{{\partial u_{\alpha \beta } \partial u_{\gamma \delta } }}u_\alpha  u_\beta  u_{\gamma \delta \xi } } \right)} _\eta \\
  &+ ( - 1)^{l + 1} \sum\limits_{\alpha ,\beta } {\left( {\frac{{\partial \sigma _{l + 2} (D^2 u)}}{{\partial u_{\alpha \beta } }}\left( {u_{\alpha \xi } u_\beta   + u_\alpha  u_{\beta \xi } } \right)} \right)} _\eta
\end{align*}

Denote by $a^{\xi \eta } = (1+|\nabla u|^2)\delta_{\xi \eta}-u_\xi u_\eta$ . Direct calculations show that
\begin{align}\label {Expression}
\begin{split}
 \sum\limits_{\xi ,\eta } {a^{\xi \eta } \varphi _{\xi \eta }  = I + II + III + IV + V} ,
 \end{split}
\end{align}
where
 \begin{align}\label {Expression-I-V}
\begin{split}
 I = & ( - 1)^{l + 1} \sum\limits_{\alpha ,\beta ,\gamma ,\delta ,\tau ,\theta ,\xi ,\eta } {a^{\xi \eta } \frac{{\partial ^3 \sigma _{l + 2} (D^2 u)}}{{\partial u_{\alpha \beta } \partial u_{\gamma \delta } \partial u_{\tau \theta } }}u_\alpha  u_\beta  } u_{\gamma \delta \xi } u_{\tau \theta \eta } ; \\
 II =& 4( - 1)^{l + 1} \sum\limits_{\alpha ,\beta ,\gamma ,\delta ,\xi ,\eta } {a^{\xi \eta } \frac{{\partial ^2 \sigma _{l + 2} (D^2 u)}}{{\partial u_{\alpha \beta } \partial u_{\gamma \delta } }}u_{\alpha \eta } u_\beta  } u_{\gamma \delta \xi } ; \\
 III =& ( - 1)^{l + 1} \sum\limits_{\alpha ,\beta ,\gamma ,\delta ,\xi ,\eta } {a^{\xi \eta } \frac{{\partial ^2 \sigma _{l + 2} (D^2 u)}}{{\partial u_{\alpha \beta } \partial u_{\gamma \delta } }}u_\alpha  u_\beta  } u_{\gamma \delta \xi \eta } ; \\
 IV =& 2( - 1)^{l + 1} \sum\limits_{\alpha ,\beta ,\xi ,\eta } {a^{\xi \eta } \frac{{\partial \sigma _{l + 2} (D^2 u)}}{{\partial u_{\alpha \beta } }}u_{\alpha \xi \eta } u_\beta  } ; \\
 V =& 2( - 1)^{l + 1} \sum\limits_{\alpha ,\beta ,\xi ,\eta } {a^{\xi \eta } \frac{{\partial \sigma _{l + 2} (D^2 u)}}{{\partial u_{\alpha \beta } }}u_{\alpha \xi } u_{\beta \eta } } .
 \end{split}
\end{align}

In the following, we come to deal with the above five terms one by one.

For $I$, we compute step by step and will arrive at \eqref{IOver}.

\begin{align}\label {I}
\begin{split}
 I =& ( - 1)^{l + 1} \sum\limits_{\alpha ,\beta ,\gamma ,\delta ,\tau ,\theta ,\xi ,\eta } {a^{\xi \eta } \frac{{\partial ^3 \sigma _{l + 2} (D^2 u)}}{{\partial u_{\alpha \beta } \partial u_{\gamma \delta } \partial u_{\tau \theta } }}u_\alpha  u_\beta  } u_{\gamma \delta \xi } u_{\tau \theta \eta }  \\
  =& ( - 1)^{l + 1} u_n ^2 \sum\limits_\xi  {a^{\xi \xi } \sum\limits_{i,j,k,l} {\frac{{\partial ^2 \sigma _{l + 1} (u_{pq} )}}{{\partial u_{ij} \partial u_{kl} }}} } u_{ij\xi } u_{kl\xi }  \\
  =& ( - 1)^{l + 1} u_n ^2 \sum\limits_\xi  {a^{\xi \xi } \sum\limits_{i \ne j} {\sigma _{l - 1} \left( {\mu \left| {ij} \right.} \right)\left( {u_{ii\xi } u_{jj\xi }  - u_{ij\xi } ^2 } \right)} }  \\
  =& u_n ^{l + 1} \sum\limits_\xi  {a^{\xi \xi } \sum\limits_{i \ne j} {\sigma _{l - 1} \left( {\lambda \left| {ij} \right.} \right)\left( {u_{ii\xi } u_{jj\xi }  - u_{ij\xi } ^2 } \right)} }  \\
  =& u_n ^{l + 1} \sum\limits_\xi  {a^{\xi \xi } \left( {\sum\limits_{i ,j \in G,i \ne j} {} {\rm{ + }}\sum\limits_{i \in G,j \in B} {}  + \sum\limits_{j \in G,i \in B} {}  + \sum\limits_{i,j \in B,i \ne j} {} } \right)\sigma _{l - 1} \left( {\lambda \left| {ij} \right.} \right)\left( {u_{ii\xi } u_{jj\xi }  - u_{ij\xi } ^2 } \right)}  \\
  = & I_1  + I_2  + I_3  + I_4 .
  \end{split}
\end{align}

It is easy to conclude that
\begin{align}\label {I1}
\begin{split}
 I_1\sim 0.
 \end{split}
\end{align}

For the term $I_2+I_3$, it follows that
\begin{align}\label {I2I3}
\begin{split}
 I_2  + I_3  =& 2u_n ^{l + 1} \sum\limits_\xi  {a^{\xi \xi } \sum\limits_{i \in G,j \in B} {\sigma _{l - 1} \left( {\lambda \left| {ij} \right.} \right)\left( {u_{ii\xi } u_{jj\xi }  - u_{ij\xi } ^2 } \right)} }  \\
  \sim & 2u_n ^{l - 1} \sum\limits_\xi  {a^{\xi \xi } \sum\limits_{i \in G,j \in B} {\sigma _{l - 1} \left( {G\left| i \right.} \right)\left( {\left( {u_n u_{ii\xi } } \right)\left( {u_n u_{jj\xi } } \right) - \left( {u_n u_{ij\xi } } \right)^2 } \right)} }  \\
  =& 2u_n ^{l - 1} \sum\limits_\xi  {a^{\xi \xi } \sum\limits_{i \in G,j \in B} {\sigma _{l - 1} \left( {G\left| i \right.} \right) {\left( { - u_n ^2 b_{ii,\xi }  + 2u_{ni} u_{i\xi }  + u_{n\xi } u_{ii} } \right)\left( { - u_n ^2 b_{jj,\xi }  + 2u_{nj} u_{j\xi } } \right)} } }  \\
  &- 2u_n ^{l - 1} \sum\limits_\xi  {a^{\xi \xi } \sum\limits_{i \in G,j \in B} {\sigma _{l - 1} \left( {G\left| i \right.} \right)} } \left( { - u_n ^2 b_{ij,\xi }  + u_{ni} u_{j\xi }  + u_{nj} u_{i\xi } } \right)^2  \\
  =& I_{21}  + I_{22}.
 \end{split}
\end{align}

By \eqref{1-ordercondition}, we have
\begin{align}\label {I21}
\begin{split}
 I_{21}  =& 4u_n ^{l - 1} \sum\limits_\xi  {a^{\xi \xi } \sum\limits_{i \in G,j \in B} {\sigma _{l - 1} \left( {G\left| i \right.} \right)\left( { - u_n ^2 b_{ii,\xi }  + 2u_{ni} u_{i\xi }  + u_{n\xi } u_{ii} } \right)u_{nj} u_{j\xi } } }  \\
  =& 4u_n ^{l - 1} \sum\limits_{i \in G,j \in B} {\sigma _{l - 1} \left( {G\left| i \right.} \right)\left( { - u_n ^2 b_{ii,n}  + 2u_{ni} ^2  + u_{nn} u_{ii} } \right)u_{nj} ^2 }  \\
   =&  - 4u_n ^{l + 1} \sum\limits_{i \in G,j \in B} {\sigma _{l - 1} \left( {G\left| i \right.} \right)b_{ii,n} } u_{nj} ^2  + 8u_n ^{l - 1} \sum\limits_{i \in G,j \in B} {\sigma _{l - 1} \left( {G\left| i \right.} \right)u_{ni} ^2 u_{nj} ^2 }\\
    &- 4lu_n ^l u_{nn} \sigma _l \left( G \right)\sum\limits_{j \in B} {u_{nj} ^2 },
 \end{split}
\end{align}
and for the second term $I_{22}$, we get
\begin{align}\label {I22}
\begin{split}
 I_{22}
  =&  - 2u_n ^{l + 3} \sum\limits_\xi  {a^{\xi \xi } \sum\limits_{i \in G,j \in B} {\sigma _{l - 1} \left( {G\left| i \right.} \right)b_{ij,\xi } ^2 } }  \\
  &+ 8u_n ^{l + 1} \sum\limits_{i \in G,j \in B} {\sigma _{l - 1} \left( {G\left| i \right.} \right)u_{ni} u_{nj} b_{ij,n} }  + 4u_n ^{l + 1} \left( {1 + u_n ^2 } \right)\sum\limits_{i \in G,j \in B} {\sigma _{l - 1} \left( {G\left| i \right.} \right)u_{nj} u_{ii} b_{ij,i} }  \\
  &- 8u_n ^{l - 1} \sum\limits_{i \in G,j \in B} {\sigma _{l - 1} \left( {G\left| i \right.} \right)u_{ni} ^2 u_{nj} ^2 }  - 2u_n ^{l - 1} \left( {1 + u_n ^2 } \right)\sum\limits_{i \in G,j \in B} {\sigma _{l - 1} \left( {G\left| i \right.} \right)u_{ii} ^2 u_{nj} ^2 }  \\
  =&  - 2u_n ^{l + 3} \sum\limits_\xi  {a^{\xi \xi } \sum\limits_{i \in G,j \in B} {\sigma _{l - 1} \left( {G\left| i \right.} \right)b_{ij,\xi } ^2 } }  \\
  &+ 8u_n ^{l + 1} \sum\limits_{i \in G,j \in B} {\sigma _{l - 1} \left( {G\left| i \right.} \right)u_{ni} u_{nj} } b_{ij,n}  - 4u_n ^{l + 2} \left( {1 + u_n ^2 } \right)\sigma _l \left( G \right)\sum\limits_{i \in G,j \in B} {u_{nj} b_{ij,i} }  \\
  &- 8u_n ^{l - 1} \sum\limits_{i \in G,j \in B} {\sigma _{l - 1} \left( {G\left| i \right.} \right)u_{ni} ^2 u_{nj} ^2 }  - 2u_n ^{l + 1} \left( {1 + u_n ^2 } \right)\sigma _l \left( G \right)\sigma _1 \left( G \right)\sum\limits_{j \in B} {u_{nj} ^2 }  \,.
 \end{split}
\end{align}

Therefore, by \eqref{I2I3},\eqref{I21} and \eqref{I22} we get

\begin{align}\label {I2+I3}
\begin{split}
 I_2  + I_3  = & - 2u_n ^{l + 3} \sum\limits_\xi  {a^{\xi \xi } \sum\limits_{i \in G,j \in B} {\sigma _{l - 1} \left( {G\left| i \right.} \right)b_{ij,\xi } ^2 } }  - 4u_n ^{l + 1} \sum\limits_{i \in G,j \in B} {\sigma _{l - 1} \left( {G\left| i \right.} \right)} u_{nj} ^2 b_{ii,n}  \\
  &+ 8u_n ^{l + 1} \sum\limits_{i \in G,j \in B} {\sigma _{l - 1} \left( {G\left| i \right.} \right)u_{ni} u_{nj} } b_{ij,n}  - 4u_n ^{l + 2} \left( {1 + u_n ^2 } \right)\sigma _l \left( G \right)\sum\limits_{i \in G,j \in B} {u_{nj} b_{ii,j} }  \\
  &- 2u_n ^{l + 1} \left( {1 + u_n ^2 } \right)\sigma _l \left( G \right)\sigma _1 \left( G \right)\sum\limits_{j \in B} {u_{nj} ^2 }  - 4lu_n ^l u_{nn} \sigma _l \left( G \right)\sum\limits_{j \in B} {u_{nj} ^2 }.
 \end{split}
\end{align}

For the term $I_4$,
\begin{align}\label {I4}
\begin{split}
 I_4  \sim& u_n ^{l - 1}\sigma _{l - 1} \left( G \right) \sum\limits_\xi  {a^{\xi \xi }\sum\limits_{i,j \in B,i \ne j} {\left( { - u_n ^2 b_{ii,\xi }  + 2u_{ni} u_{i\xi } } \right)\left( { - u_n ^2 b_{jj,\xi }  + 2u_{nj} u_{j\xi } } \right)} }  \\
  &- u_n ^{l - 1}\sigma _{l - 1} \left( G \right) \sum\limits_\xi  {a^{\xi \xi }\sum\limits_{i,j \in B,i \ne j} {\left( { - u_n ^2 b_{ij,\xi }  + u_{ni} u_{j\xi }  + u_{nj} u_{i\xi } } \right)^2 } }  \\
  =& I_{41}  + I_{42} .
 \end{split}
\end{align}

According to \eqref{1-ordercondition} we derive
\begin{align}\label {I41}
\begin{split}
 I_{41}  =& u_n ^{l - 1} \sigma _{l - 1} \left( G \right)\sum\limits_\xi  {a^{\xi \xi }\sum\limits_{i,j \in B,i \ne j} {\left( { - u_n ^2 b_{ii,\xi }  + 2u_{ni} u_{i\xi } } \right)\left( { - u_n ^2 b_{jj,\xi }  + 2u_{nj} u_{j\xi } } \right)} }  \\
  =& u_n ^{l - 1} \sigma _{l - 1} \left( G \right)\sum\limits_\xi  {a^{\xi \xi } \sum\limits_{i,j \in B,i \ne j} {\left[ {u_n ^4 b_{ii,\xi } b_{jj,\xi }  - 4u_n ^2 u_{nj} u_{j\xi } b_{ii,\xi }  + 4u_{ni} u_{nj} u_{i\xi } u_{j\xi } } \right]} }  \\
  = & - u_n ^{l + 3} \sigma _{l - 1} \left( G \right)\sum\limits_{i \in B} {b_{ii,n} ^2 }  - u_n ^{l + 3} \left( {1 + u_n ^2 } \right)\sigma _{l - 1} \left( G \right)\sum\limits_{j = 1}^{n - 1} {\sum\limits_{i \in B} {b_{ii,j} ^2 } }  \\
  &+ 4u_n ^{l + 1} \sigma _{l - 1} \left( G \right)\sum\limits_{j \in B} {u_{nj} ^2 b_{jj,n} }  + 4u_n ^{l - 1} \sigma _{l - 1} \left( G \right)\sum\limits_{i,j \in B,i \ne j} {u_{ni} ^2 u_{nj} ^2 }
 \end{split}
\end{align}
and
\begin{align}\label {I42}
\begin{split}
 I_{42}  =&  - u_n ^{l - 1} \sum\limits_\xi  {a^{\xi \xi } \sigma _{l - 1} \left( G \right)\sum\limits_{i,j \in B,i \ne j} {\left( { - u_n ^2 b_{ij,\xi }  + u_{ni} u_{j\xi }  + u_{nj} u_{i\xi } } \right)^2 } }  \\
  =&  - u_n ^{l - 1} \sigma _{l - 1} \left( G \right)\sum\limits_{i,j \in B,i \ne j} {\left( { - u_n ^2 b_{ij,n}  + 2u_{ni} u_{nj} } \right)^2 }  - u_n ^{l + 3} \left( {1 + u_n ^2 } \right)\sigma _{l - 1} \left( G \right)\sum\limits_{k = 1}^{n - 1} {\sum\limits_{i,j \in B,i \ne j} {b_{ij,k} ^2 } }  \\
  =&  - u_n ^{l + 3} \sigma _{l - 1} \left( G \right)\sum\limits_{i,j \in B,i \ne j} {b_{ij,n} ^2 }  + 4u_n ^{l + 1} \sigma _{l - 1} \left( G \right)\sum\limits_{i,j \in B,i \ne j} {u_{ni} u_{nj} b_{ij,n} }  \\
  &- 4u_n ^{l - 1} \sigma _{l - 1} \left( G \right)\sum\limits_{i,j \in B,i \ne j} {u_{ni} ^2 u_{nj} ^2 }  - u_n ^{l + 3} \left( {1 + u_n ^2 } \right)\sigma _{l - 1} \left( G \right)\sum\limits_{k = 1}^{n - 1} {\sum\limits_{i,j \in B,i \ne j} {b_{ij,k} ^2 } }\ .
  \end{split}
\end{align}

Then we can conclude from \eqref{I4}, \eqref{I41} and \eqref{I42} that
\begin{align}\label {I4Over}
\begin{split}
  I_4  = & - u_n ^{l + 3} \sigma _{l - 1} \left( G \right)\sum\limits_{i,j \in B} {b_{ij,n} ^2 }  \\
  &- u_n ^{l + 3} \left( {1 + u_n ^2 } \right)\sigma _{l - 1} \left( G \right)\sum\limits_{k = 1}^{n - 1} {\sum\limits_{i,j \in B} {b_{ij,k} ^2 } }  + 4u_n ^{l + 1} \sigma _{l - 1} \left( G \right)\sum\limits_{i,j \in B} {u_{ni} u_{nj} b_{ij,n} } \ .
\end{split}
\end{align}

 Now, together with \eqref{I}, \eqref{I1}, \eqref{I2+I3} and \eqref{I4Over} we have
 \begin{align}\label {IOver}
\begin{split}
 I \sim & - 2u_n ^{l + 3} \sum\limits_\xi  {a^{\xi \xi } \sum\limits_{i \in G,j \in B} {\sigma _{l - 1} \left( {G\left| i \right.} \right)b_{ij,\xi } ^2 } }  - 4u_n ^{l + 1} \sum\limits_{i \in G,j \in B} {\sigma _{l - 1} \left( {G\left| i \right.} \right)} u_{nj} ^2 b_{ii,n}  \\
  &+ 8u_n ^{l + 1} \sum\limits_{i \in G,j \in B} {\sigma _{l - 1} \left( {G\left| i \right.} \right)u_{ni} u_{nj} } b_{ij,n}  - 4u_n ^{l + 2} \left( {1 + u_n ^2 } \right)\sigma _l \left( G \right)\sum\limits_{i \in G,j \in B} {u_{nj} b_{ii,j} }  \\
  &- 2u_n ^{l + 1} \left( {1 + u_n ^2 } \right)\sigma _l \left( G \right)\sigma _1 \left( G \right)\sum\limits_{j \in B} {u_{nj} ^2 }  - 4lu_n ^l u_{nn} \sigma _l \left( G \right)\sum\limits_{j \in B} {u_{nj} ^2 }  \\
  &- u_n ^{l + 3} \sigma _{l - 1} \left( G \right)\sum\limits_{i,j \in B} {b_{ij,n} ^2 }  - u_n ^{l + 3} \left( {1 + u_n ^2 } \right)\sigma _{l - 1} \left( G \right)\sum\limits_{k = 1}^{n - 1} {\sum\limits_{i,j \in B} {b_{ij,k} ^2 } } \\
   &+ 4u_n ^{l + 1} \sigma _{l - 1} \left( G \right)\sum\limits_{i,j \in B} {u_{ni} u_{nj} b_{ij,n} }\ .
\end{split}
\end{align}

For the term $II$, we will arrive at \eqref{IIOver}.
\begin{align}\label {II}
\begin{split}
 II = &4( - 1)^{l + 1} \sum\limits_\xi  {a^{\xi \xi } } \sum\limits_{\alpha, \beta, \gamma, \delta } {\frac{{\partial ^2 \sigma _{l + 2} (D^2 u)}}{{\partial u_{\alpha \beta } \partial u_{\gamma \delta } }}u_{\alpha \xi } u_\beta  u_{\gamma \delta \xi } }  \\
  = &4( - 1)^{l + 1} u_n \sum\limits_\xi  {a^{\xi \xi } } \sum\limits_{\gamma ,\delta } {\frac{{\partial ^2 \sigma _{l + 2} (D^2 u)}}{{\partial u_{nn} \partial u_{\gamma \delta } }}u_{n\xi } u_{\gamma \delta \xi } } \\
   &+ 4( - 1)^{l + 1} u_n \sum\limits_\xi  {a^{\xi \xi } } \sum\limits_i {\sum\limits_{\gamma ,\delta } {\frac{{\partial ^2 \sigma _{l + 2} (D^2 u)}}{{\partial u_{in} \partial u_{\gamma \delta } }}u_{i\xi } u_{\gamma \delta \xi } } }
  = II_1  + II_2.
\end{split}
\end{align}

For the first term $II_1$ we have by \eqref{bij-xi},
\begin{align}\label {II1}
\begin{split}
 II_1  = &4( - 1)^{l + 1} u_n \sum\limits_\xi  {a^{\xi \xi } } \sum\limits_{\gamma ,\delta } {\frac{{\partial ^2 \sigma _{l + 2} (D^2 u)}}{{\partial u_{nn} \partial u_{\gamma \delta } }}u_{n\xi } u_{\gamma \delta \xi } }  \\
  =& 4( - 1)^{l + 1} u_n \sum\limits_\xi  {a^{\xi \xi } } \sum\limits_i {\sigma _l \left( {\mu \left| i \right.} \right)u_{n\xi } u_{ii\xi } }  \\
 \sim & - 4u_n ^{l + 1} \sum\limits_\xi  {a^{\xi \xi } } \sum\limits_{i \in B} {\sigma _l \left( G \right)} u_{n\xi } u_{ii\xi }  \\
  =&  - 4u_n ^l \sigma _l \left( G \right)\sum\limits_\xi  {a^{\xi \xi } } \sum\limits_{i \in B} {u_{n\xi } ( - u_n ^2 b_{ii,\xi }  + 2u_{ni} u_{i\xi } )}  \\
  = & - 8fu_n ^l \sigma _l \left( G \right)\sum\limits_{j \in B} {u_{nj} ^2 }  - 8u_n ^{l + 1} \left( {1 + u_n ^2 } \right)\sigma _l \left( G \right)\sigma _1 \left( G \right)\sum\limits_{j \in B} {u_{nj} ^2 }.
 \end{split}
\end{align}

In the following, we come to settle the second term $II_2$.
\begin{align}\label {II2}
\begin{split}
II_2 = & 4( - 1)^{l + 1} u_n \sum\limits_\xi  {a^{\xi \xi } } \frac{{\partial ^2 \sigma _{l + 2} (D^2 u)}}{{\partial u_{in} \partial u_{\gamma \delta } }}u_{i\xi } u_{\gamma \delta \xi }  \\
  =& 4( - 1)^{l + 1} u_n \sum\limits_\xi  {a^{\xi \xi } } \sum\limits_{i = 1}^{n - 1} {\frac{\partial }{{\partial u_{\gamma \delta } }}\left( { - u_{ni} \sigma _l (u_{pq} \left| i \right.) + \sum\limits_{j \ne i} {u_{nj} u_{ij} } \sigma _{l - 1} (u_{pq} \left| {ij} \right.)} \right)} u_{i\xi } u_{\gamma \delta \xi }  \\
  =& II_{21}  + II_{22} .
\end{split}
\end{align}

It follows that
\begin{align}\label {II21}
\begin{split}
II_{21}  =& 4( - 1)^{l + 1} u_n \sum\limits_\xi  {a^{\xi \xi } } \sum\limits_{i = 1}^{n - 1} {\frac{\partial }{{\partial u_{\gamma \delta } }}\left( { - u_{ni} \sigma _l (u_{pq} \left| i \right.)} \right)} u_{i\xi } u_{\gamma \delta \xi }  \\
  =& 4( - 1)^{l + 1} u_n \sum\limits_\xi  {a^{\xi \xi } } \sum\limits_{i = 1}^{n - 1} {\left( { - \sigma _l (u_{pq} \left| i \right.)} \right)} u_{i\xi } u_{ni\xi }  \\
  &- 4( - 1)^{l + 1} u_n \sum\limits_\xi  {a^{\xi \xi } } \sum\limits_{i \ne j} {\sigma _{l - 1} (u_{pq} \left| i j\right.)u_{ni} } u_{i\xi } u_{jj\xi }
  = II_{211}  + II_{212} ;
\end{split}
\end{align}
For $II_{211}$, we derive that
\begin{align}\label {II211}
\begin{split}
 II_{211}  =& 4( - 1)^{l + 1} u_n \sum\limits_\xi  {a^{\xi \xi } } \sum\limits_{i = 1}^{n - 1} {\left( { - \sigma _l (u_{pq} \left| i \right.)} \right)} u_{i\xi } u_{ni\xi }  \\
  \sim &4u_n ^{l + 1} \sigma _l \left( G \right)\sum\limits_\xi  {a^{\xi \xi } } \sum\limits_{i \in B} {u_{i\xi } u_{ni\xi } }  \\
  \sim& 4u_n ^{l + 1} \sigma _l \left( G \right)\sum\limits_{i \in B} {u_{in} u_{nin} }
  =4u_n ^{l + 1} \sigma _l \left( G \right)\sum\limits_{i \in B} {u_{in} u_{nni} }\ .
\end{split}
\end{align}

We now pause for a while to derive a observation for the proceed of the whole proof.

Differentiate both sides of the equation $ a^{\alpha\beta}u_{\alpha\beta} =f$ with respect to $e_i(i\in B)$,
\begin{align*}
\sum\limits_{\alpha,\beta}{a^{\alpha \beta }} _{,i} u_{\alpha \beta }  +\sum\limits_{\alpha,\beta} a^{\alpha \beta } u_{\alpha \beta i}  = f_i.
\end{align*}

Easy to get that for $i\in B$,
\begin{align}\label {daoshuzhuanyi}
\begin{split}
\sum\limits_{\alpha,\beta}{a^{\alpha \beta }} _{,i} u_{\alpha \beta }=2u_nu_{ni}\triangle u-2u_n\sum\limits_{\alpha}u_{n\alpha}u_{i\alpha}=-2u_n^2u_{ni}\sigma_1(G).
\end{split}
\end{align}

Therefore, by \eqref{bij-xi} the term $u_{nni}$ now becomes
\begin{align}\label {unni}
\begin{split}
u_{nni}= f_i +2u_n^2u_{ni}\sigma_1(G) - u_n ^{ - 1} \left( {1 + u_n ^2 } \right)\sum\limits_{j \in G} {\left( { - u_n ^2 b_{jj,i}  + u_{ni} u_{jj} } \right)}.
\end{split}
\end{align}

Plug this expression into \eqref{II211}, we have
\begin{align}\label {II211Over}
\begin{split}
II_{211}   =& 4u_n ^{l + 1} \sigma _l \left( G \right)\sum\limits_{j \in B} {u_{nj} f_j } + 4u_n ^{l + 1} \sigma _l \left( G \right)\sigma _1 \left( G \right)\left( {1 + 3u_n ^2 } \right)\sum\limits_{j \in B} {u_{nj} ^2 }  \\
 &+ 4u_n ^{l + 2} \sigma _l \left( G \right)\left( {1 + u_n ^2 } \right)\sum\limits_{i \in G,j \in B} {u_{nj} b_{ii,j} }.
\end{split}
\end{align}

For the term $II_{212}$,
\begin{align}\label {}
\begin{split}
II_{212}  = & - 4( - 1)^{l + 1} u_n \sum\limits_\xi  {a^{\xi \xi } } \sum\limits_{i \ne j} {\sigma _{l - 1} (u_{pq} \left| ij\right.)u_{ni} } u_{i\xi } u_{jj\xi }  \\
  =& 4( - 1)^{l} u_n \sum\limits_\xi  {a^{\xi \xi } } \left( {\sum\limits_{i,j \in G,i \ne j} {} {\rm{ + }}\sum\limits_{i \in G,j \in B} {}  + \sum\limits_{j \in G,i \in B} {}  + \sum\limits_{i,j \in B,i \ne j} {} } \right)\sigma _{l - 1} (\mu \left| ij\right.)u_{ni} u_{i\xi } u_{jj\xi } \\
  =& II_{2121}  + II_{2122}  + II_{2123}  + II_{2124} .
\end{split}
\end{align}

 Obvious to find that
 \begin{align}\label {II2121}
\begin{split}
II_{2121}\sim0.
\end{split}
\end{align}

By \eqref{bij-xi} and \eqref{1-ordercondition}, we get
\begin{align}\label {II2122}
\begin{split}
 II_{2122}  =& 4( - 1)^l \sum\limits_\xi  {a^{\xi \xi } } \sum\limits_{i \in G,j \in B} {\sigma _{l - 1} (\mu \left| ij\right.)u_{ni} u_{i\xi } \left( { - u_n ^2 b_{jj,\xi }  + 2u_{nj} u_{j\xi }  + u_{n\xi } u_{jj} } \right)}  \\
 \sim& - 4u_n ^{l - 1} \sum\limits_\xi  {a^{\xi \xi } } \sum\limits_{i \in G,j \in B} {\sigma _{l - 1} (\lambda \left| i \right.)u_{ni} u_{i\xi } \left( { - u_n ^2 b_{jj,\xi }  + 2u_{nj} u_{j\xi } } \right)}  \\
=& - 8u_n ^{l - 1} \sum\limits_{i \in G,j \in B} {\sigma _{l - 1} (G\left| i \right.)u_{ni} ^2 u_{nj} ^2 } ,
\end{split}
\end{align}
similarly,
\begin{align}\label {II2123}
\begin{split}
 II_{2123}  = &4( - 1)^l \sum\limits_\xi  {a^{\xi \xi } } \sum\limits_{j \in G,i \in B} {\sigma _{l - 1} (\mu \left| {ij} \right.)u_{ni} u_{i\xi } \left( { - u_n ^2 b_{jj,\xi }  + 2u_{nj} u_{j\xi }  + u_{n\xi } u_{jj} } \right)}  \\
  \sim& 4u_n ^{l + 1} \sum\limits_{i \in G,j \in B} {\sigma _{l - 1} (\lambda \left| i \right.)u_{nj} ^2 b_{ii,n} }  - 8u_n ^{l - 1} \sum\limits_{i \in G,j \in B} {\sigma _{l - 1} (\lambda \left| i \right.)u_{ni} ^2 u_{nj} } ^2  \\
  &+ 4lu_n ^l u_{nn} \sigma _l \left( G \right)\sum\limits_{j \in B} {u_{nj} } ^2 \ ,
\end{split}
\end{align}
and
\begin{align}\label {II2124}
\begin{split}
 II_{2124}  =& 4( - 1)^l \sum\limits_\xi  {a^{\xi \xi } } \sum\limits_{i,j \in B,i \ne j} {\sigma _{l - 1} (\mu \left| ij\right.)u_{ni} u_{i\xi } \left( { - u_n ^2 b_{jj,\xi }  + 2u_{nj} u_{j\xi }  + u_{n\xi } u_{jj} } \right)}  \\
  = & - 4u_n ^{l - 1} \sum\limits_{i,j \in B,i \ne j} {\sigma _{l - 1} (\lambda \left| ij\right.)u_{in} ^2 \left( { - u_n ^2 b_{jj,n}  + 2u_{nj} ^2 } \right)}  \\
  \sim &4u_n ^{l + 1} \sigma _{l - 1} (G)\sum\limits_{i,j \in B,i \ne j} {u_{in} ^2 b_{jj,n} }  - 8u_n ^{l - 1} \sigma _{l - 1} (G)\sum\limits_{i,j \in B,i \ne j} {u_{in} ^2 u_{nj} ^2 }  \\
  \sim & - 4u_n ^{l + 1} \sigma _{l - 1} (G)\sum\limits_{j \in B} {u_{nj} ^2 b_{jj,n} }  - 8u_n ^{l - 1} \sigma _{l - 1} (G)\sum\limits_{i,j \in B,i \ne j} {u_{ni} ^2 u_{nj} ^2 }\ .
 \end{split}
\end{align}

Therefore, according to \eqref{II2121}, \eqref{II2122}, \eqref{II2123} and \eqref{II2124} we have
\begin{align}\label {II212Over}
\begin{split}
 II_{212} \sim &4u_n ^{l + 1} \sum\limits_{i \in G,j \in B} {\sigma _{l - 1} (G\left| i \right.)u_{nj} ^2 b_{ii,n} }  - 4u_n ^{l + 1} \sigma _{l - 1} (G)\sum\limits_{j \in B} {u_{nj} ^2 b_{jj,n} }  \\
  &- 16u_n ^{l - 1} \sum\limits_{i \in G,j \in B} {\sigma _{l - 1} (G\left| i \right.)u_{ni} ^2 u_{nj} } ^2  - 8u_n ^{l - 1} \sigma _{l - 1} (G)\sum\limits_{i,j \in B,i \ne j} {u_{ni} ^2 u_{nj} ^2 }  \\
  & + 4lu_n ^l u_{nn} \sigma _l \left( G \right)\sum\limits_{j \in B} {u_{nj} } ^2  \ .
\end{split}
\end{align}

By \eqref{II211Over} and \eqref{II212Over}, we get that
\begin{align}\label {II21Over}
\begin{split}
 II_{21} \sim&4u_n ^{l + 1} \sigma _l \left( G \right)\sum\limits_{j \in B} {u_{nj} f_j }  + 4u_n ^{l + 1} \sigma _l \left( G \right)\sigma _1 \left( G \right)\left( {1 + 3u_n ^2 } \right)\sum\limits_{j \in B} {u_{nj} ^2 }  \\
  &+ 4u_n ^{l + 2} \sigma _l \left( G \right)\left( {1 + u_n ^2 } \right)\sum\limits_{i \in G,j \in B} {u_{nj} b_{ii,j} }  + 4u_n ^{l + 1} \sum\limits_{i \in G,j \in B} {\sigma _{l - 1} (G\left| i \right.)u_{nj} ^2 b_{ii,n} }  \\
  &- 4u_n ^{l + 1} \sigma _{l - 1} (G)\sum\limits_{j \in B} {u_{nj} ^2 b_{jj,n} }  - 16u_n ^{l - 1} \sum\limits_{i \in G,j \in B} {\sigma _{l - 1} (G\left| i \right.)u_{ni} ^2 u_{nj} } ^2  \\
  &+ 4lu_n ^l u_{nn} \sigma _l \left( G \right)\sum\limits_{j \in B} {u_{nj} } ^2  - 8u_n ^{l - 1} \sigma _{l - 1} (G)\sum\limits_{i,j \in B,i\neq,j} {u_{ni} ^2 u_{nj} ^2 }  \ .
\end{split}
\end{align}

In the sequel, we come to settle the term $II_{22}$.
\begin{align}\label {II22}
\begin{split}
II_{22}  =& 4\left( { - 1} \right)^{l + 1} u_n \sum\limits_{\xi  = 1}^n {a^{\xi \xi } } \sum\limits_{i \ne j} {\sum\limits_{\gamma, \delta } {\frac{\partial }{{\partial u_{\gamma \delta } }}} } \left[ {u_{nj} u_{ji} \sigma _{l - 1} \left( {u_{pq} \left| {ij} \right.} \right)} \right]u_{i\xi } u_{\gamma\delta\xi }  \\
  = &4u_n ^l \sum\limits_{\xi  = 1}^n {a^{\xi \xi } } \left( {\sum\limits_{i,j \in G,i \ne j} {} {\rm{ + }}\sum\limits_{i \in G,j \in B} {}  + \sum\limits_{j \in G,i \in B} {}  + \sum\limits_{i,j \in B,i \ne j} {} } \right)u_{nj} \sigma _{l - 1} \left( {\lambda \left| {ij} \right.} \right)u_{i\xi } u_{ji\xi } \\
  =&II_{221}+II_{222}+II_{223}+II_{224} .
\end{split}
\end{align}

Step by step, we will give the final expression of $II_{22}$. Firstly,
\begin{align}\label {II221}
\begin{split}
II_{221} \sim 0.
\end{split}
\end{align}

Secondly, by \eqref{bij-xi} and Lemma \ref{Derivative comutation},
\begin{align}\label {II222}
\begin{split}
 II_{222} \sim & 4u_n ^l \sum\limits_{\xi  = 1}^n {a^{\xi \xi } } \sum\limits_{i \in G,j \in B} {u_{nj} \sigma _{l - 1} \left( {G\left| i \right.} \right)u_{i\xi } u_{ij\xi } }  \\
  =& 4u_n ^l \sum\limits_{i \in G,j \in B} {\sigma _{l - 1} \left( {G\left| i \right.} \right)u_{nj} u_{ni} u_{ijn} }  + 4u_n ^l \left( {1 + u_n ^2 } \right)\sum\limits_{i \in G,j \in B} {\sigma _{l - 1} \left( {G\left| i \right.} \right)u_{nj} u_{ii} u_{iji} }  \\
=&4u_n ^l \sum\limits_{i \in G,j \in B} {\sigma _{l - 1} \left( {G\left| i \right.} \right)u_{nj} u_{ni} u_{ijn} }  - 4u_n ^l \left( {1 + u_n ^2 } \right)\sum\limits_{i \in G,j \in B} {\sigma _l \left( G \right)u_{nj} \left[ { - u_n ^2 b_{ii,j}  + u_{nj} u_{ii} } \right]}  \\
  = &4u_n ^l \sum\limits_{i \in G,j \in B} {\sigma _{l - 1} \left( {G\left| i \right.} \right)u_{nj} u_{ni} u_{ijn} }  + 4u_n ^{l + 2} \left( {1 + u_n ^2 } \right)\sigma _l \left( G \right)\sum\limits_{i \in G,j \in B} {u_{nj} b_{ii,j} }  \\
  &+ 4u_n ^{l + 1} \left( {1 + u_n ^2 } \right)\sigma _l \left( G \right)\sigma _1 \left( G \right)\sum\limits_{j \in B} {u_{nj} ^2 }\ .  \\
\end{split}
\end{align}

A similar process shows that
\begin{align}\label {II223}
\begin{split}
II_{223} \sim 4u_n ^l \sum\limits_{i \in G,j \in B} {\sigma _{l - 1} \left( {G\left| i \right.} \right)u_{nj} u_{ni} u_{ijn} }\ ,
\end{split}
\end{align}
and
\begin{align}\label {II224}
\begin{split}
II_{224}  \sim 4u_n ^l \sigma _{l - 1} \left( G \right)\sum\limits_{i,j \in B,i \ne j} {u_{nj} u_{ni} u_{ijn} }\ .
\end{split}
\end{align}

Combining \eqref{II221}$-$\eqref{II224}, we have
\begin{align}\label {II22heqilai}
\begin{split}
 II_{22} \sim& 8 u_n ^l \sum\limits_{i \in G,j \in B} {\sigma _{l - 1} \left( {G\left| i \right.} \right)u_{nj} u_{ni} u_{ijn} } {\rm{ + }}4u_n ^l \sigma _{l - 1} \left( G \right)\sum\limits_{i,j \in B,i \ne j} {u_{nj} u_{ni} u_{ijn} }  \\
  &+ 4u_n ^{l + 2} \left( {1 + u_n ^2 } \right)\sigma _l \left( G \right)\sum\limits_{i \in G,j \in B} {u_{nj} b_{ii,j} }  + 4u_n ^{l + 1} \left( {1 + u_n ^2 } \right)\sigma _l \left( G \right)\sigma _1 \left( G \right)\sum\limits_{j \in B} {u_{nj} ^2 } .
 \end{split}
\end{align}

Once again by \eqref{bij-xi}, we revise the form of $II_{22}$ into the following
\begin{align}\label {II22Over}
\begin{split}
 II_{22}  \sim & - 8 u_n ^{l + 1} \sum\limits_{i \in G,j \in B} {\sigma _{l - 1} \left( {G\left| i \right.} \right)u_{nj} u_{ni} b_{ij,n} }  - 4u_n ^{l + 1} \sigma _{l - 1} \left( G \right)\sum\limits_{i,j \in B,i \ne j} {u_{nj} u_{ni} b_{ij,n} }  \\
  &+ 16u_n ^{l - 1} \sum\limits_{i \in G,j \in B} {\sigma _{l - 1} \left( {G\left| i \right.} \right)u_{ni} ^2 u_{nj} ^2 } {\rm{ + }}8u_n ^{l - 1} \sigma _{l - 1} \left( G \right)\sum\limits_{i,j \in B,i \ne j} {u_{ni} ^2 u_{nj} ^2 }  \\
  &+ 4u_n ^{l + {\rm{2}}} \left( {1 + u_n ^2 } \right)\sigma _l \left( G \right)\sum\limits_{i \in G,j \in B} {u_{nj} b_{ii,j} }  + 4u_n ^{l + 1} \left( {1 + u_n ^2 } \right)\sigma _l \left( G \right)\sigma _1 \left( G \right)\sum\limits_{j \in B} {u_{nj} ^2 } \ .
\end{split}
\end{align}

Combining \eqref{II2}, \eqref{II21Over} and \eqref{II22Over}, we then have
\begin{align}\label {II2Over}
\begin{split}
 II_2 \sim& 4u_n ^{l + 1} \sigma _l \left( G \right)\sum\limits_{j \in B} {u_{nj} f_j }  + 8u_n ^{l + 1} \sigma _l \left( G \right)\sigma _1 \left( G \right)\left( {1 + 2u_n ^2 } \right)\sum\limits_{j \in B} {u_{nj} ^2 }  \\
  &+ 8u_n ^{l + 2} \sigma _l \left( G \right)\left( {1 + u_n ^2 } \right)\sum\limits_{i \in G,j \in B} {u_{nj} b_{ii,j} }  + 4u_n ^{l + 1} \sum\limits_{i \in G,j \in B} {\sigma _{l - 1} (G\left| i \right.)u_{nj} ^2 b_{ii,n} }  \\
  &+ 4lu_n ^l u_{nn} \sigma _l \left( G \right)\sum\limits_{j \in B} {u_{nj} } ^2  - {\rm{8}}u_n ^{l + 1} \sum\limits_{i \in G,j \in B} {\sigma _{l - 1} \left( {G\left| i \right.} \right)u_{nj} u_{ni} b_{ij,n} }  \\
  &- 4u_n ^{l + 1} \sigma _{l - 1} \left( G \right)\sum\limits_{i,j \in B} {u_{nj} u_{ni} b_{ij,n} }  \ .
\end{split}
\end{align}

Therefore, by \eqref{II}, \eqref{II1} and \eqref{II2Over},
\begin{align}\label {IIOver}
\begin{split}
II \sim&  - 8fu_n ^l \sigma _l \left( G \right)\sum\limits_{j \in B} {u_{nj} ^2 }  + 4u_n ^{l + 1} \sigma _l \left( G \right)\sum\limits_{j \in B} {u_{nj} f_j }  + 8u_n ^{l + 3} \sigma _l \left( G \right)\sigma _1 \left( G \right)\sum\limits_{j \in B} {u_{nj} ^2 }  \\
  &+ 8u_n ^{l + 2} \sigma _l \left( G \right)\left( {1 + u_n ^2 } \right)\sum\limits_{i \in G,j \in B} {u_{nj} b_{ii,j} }  + 4u_n ^{l + 1} \sum\limits_{i \in G,j \in B} {\sigma _{l - 1} (G\left| i \right.)u_{nj} ^2 b_{ii,n} }  \\
  &+ 4lu_n ^l u_{nn} \sigma _l \left( G \right)\sum\limits_{j \in B} {u_{nj} } ^2  - {\rm{8}}u_n ^{l + 1} \sum\limits_{i \in G,j \in B} {\sigma _{l - 1} \left( {G\left| i \right.} \right)u_{nj} u_{ni} b_{ij,n} }\\
  &  - 4u_n ^{l + 1} \sigma _{l - 1} \left( G \right)\sum\limits_{i,j \in B} {u_{nj} u_{ni} b_{ij,n} }\  .
\end{split}
\end{align}

So, by \eqref{IOver} and \eqref{IIOver}, it follows that
 \begin{align}\label {I+II}
\begin{split}
I + II \sim&  - 2u_n ^{l + 3} \sum\limits_\xi  {a^{\xi \xi } \sum\limits_{i \in G,j \in B} {\sigma _{l - 1} \left( {G\left| i \right.} \right)b_{ij,\xi } ^2 } }  - u_n ^{l + 3} \sigma _{l - 1} \left( G \right)\sum\limits_{i,j \in B} {b_{ij,n} ^2 }  \\
  &- u_n ^{l + 3} \left( {1 + u_n ^2 } \right)\sigma _{l - 1} \left( G \right)\sum\limits_{k = 1}^{n - 1} {\sum\limits_{i,j \in B} {b_{ij,k} ^2 } }  - 8fu_n ^l \sigma _l \left( G \right)\sum\limits_{j \in B} {u_{nj} ^2 }  \\
  &- 2u_n ^{l + 1} \left( {1 - 3u_n ^2 } \right)\sigma _l \left( G \right)\sigma _1 \left( G \right)\sum\limits_{j \in B} {u_{nj} ^2 }  + 4u_n ^{l + 1} \sigma _l \left( G \right)\sum\limits_{j \in B} {u_{jn} f_j }  \\
  &+ 4u_n ^{l + 2} \left( {1 + u_n ^2 } \right)\sigma _l \left( G \right)\sum\limits_{i \in G,j \in B} {u_{nj} b_{ii,j} } \ .
\end{split}
\end{align}

Now, it's time for us to compute the third term $III$.

By Lemma \ref{Derivative comutation}, we then have
\begin{align}\label {III}
\begin{split}
 III =& ( - 1)^{l + 1} \sum\limits_{\alpha ,\beta ,\gamma ,\delta ,\xi ,\eta } {a^{\xi \eta } \frac{{\partial ^2 \sigma _{l + 2} (D^2 u)}}{{\partial u_{\alpha \beta } \partial u_{\gamma \delta } }}u_\alpha  u_\beta  } u_{\gamma \delta \xi \eta }  \\
    =& ( - 1)^{l + 1} u_n ^2 \sum\limits_{\alpha ,\beta ,\gamma ,\delta } {a^{\alpha \beta } \frac{{\partial \sigma _{l + 1} (u_{pq} )}}{{\partial u_{\gamma \delta } }}u_{\gamma \delta \alpha \beta } }  \\
  \sim & - u_n ^{l + 2} \sigma _l \left( G \right)\sum\limits_{\alpha ,\beta } {\sum\limits_{i \in B} {a^{\alpha \beta } u_{ii\alpha \beta } } }  \\
  = & - u_n ^{l + 2} \sigma _l \left( G \right)\sum\limits_{\alpha ,\beta } {\sum\limits_{i \in B} {a^{\alpha \beta } \left( {u_{\alpha \beta ii}  + 2\sum\limits_{\xi }u_{\xi i} R_{\xi \alpha i\beta }  - 2\sum\limits_{\xi }u_{\alpha \xi } R_{\xi i\beta i} } \right)} }  \\
  = &III_1  + III_2  + III_3  \ .
\end{split}
\end{align}

For the last two terms $III_2$ and $III_3$, it is obvious to deduce that
\begin{align}\label {III2+III3}
\begin{split}
III_2\sim0,\ \ III_3\sim2\left( {n - l-1} \right)\epsilon fu_n ^{l + 2} \sigma _l \left( G \right).
\end{split}
\end{align}
Remark that $\epsilon$ is just the sectional curvature of the space form.

In the following, we come to compute the more complicated term $III_1$.
\begin{align}\label {III1}
\begin{split}
III_1  =&  - u_n ^{l + 2} \sigma _l \left( G \right)\sum\limits_{\alpha ,\beta } {\sum\limits_{i \in B} {a^{\alpha \beta } u_{\alpha \beta ii} } }  \\
  = & - u_n ^{l + 2} \sigma _l \left( G \right)\sum\limits_{\alpha ,\beta } {\sum\limits_{i \in B} {\left[ {\left( {a^{\alpha \beta } u_{\alpha \beta } } \right)_{ii}  - \left( {a^{\alpha \beta } } \right)_{ii} u_{\alpha \beta }  - 2\left( {a^{\alpha \beta } } \right)_i u_{\alpha \beta i} } \right]} }  \\
  = & - u_n ^{l + 2} \sigma _l \left( G \right)\sum\limits_{j \in B} {f_{jj} }  + III_{11}  + III_{12}  \ .
  \end{split}
\end{align}

For $III_{11}$, we have
\begin{align}\label {III11}
\begin{split}
 III_{11}  =& u_n ^{l + 2} \sigma _l \left( G \right)\sum\limits_{\alpha ,\beta } {\sum\limits_{j \in B} {\left( {a^{\alpha \beta } } \right)_{jj} u_{\alpha \beta } } }  \\
  =& u_n ^{l + 2} \sigma _l \left( G \right)\sum\limits_{\alpha ,\beta } {\sum\limits_{j \in B} {\left( {2\sum\limits_\xi  {u_\xi  u_{\xi j} \delta _{\alpha \beta } }  - 2u_{\alpha j} u_\beta  } \right)_j u_{\alpha \beta } } }  \\
  =& 2u_n ^{l + 2} \sigma _l \left( G \right)\sum\limits_{j \in B} {\sum\limits_\alpha  {u_{\alpha j} ^2 \Delta u} }  - 2u_n ^{l + 2} \sigma _l \left( G \right)\sum\limits_{\alpha ,\beta } {\sum\limits_{j \in B} {u_{\alpha j} u_{\beta j} u_{\alpha \beta } } }  \\
  &+ 2u_n ^{l + 3} \sigma _l \left( G \right)\sum\limits_{j \in B} {u_{njj} \Delta u}  - 2u_n ^{l + 3} \sigma _l \left( G \right)\sum\limits_{\alpha  } {\sum\limits_{j \in B} {u_{\alpha jj} u_{n\alpha } } }  \\
  =& III_{111}  + III_{112}  \ ,
 \end{split}
\end{align}
where,
\begin{align}\label {III111}
\begin{split}
 III_{111}  =& 2u_n ^{l + 2} \sigma _l \left( G \right)\sum\limits_{j \in B} {\sum\limits_\alpha  {u_{\alpha j} ^2 \Delta u} }  - 2u_n ^{l + 2} \sigma _l \left( G \right)\sum\limits_{\alpha ,\beta } {\sum\limits_{j \in B} {u_{\alpha j} u_{\beta j} u_{\alpha \beta } } }  \\
  =& 2u_n ^{l + 2} \sigma _l \left( G \right)\sum\limits_{j \in B} {u_{nj} ^2 \left( {\Delta u - u_{nn} } \right)}  =  - 2u_n ^{l + 3} \sigma _l \left( G \right)\sigma _1 \left( G \right)\sum\limits_{j \in B} {u_{nj} ^2 }  \\
 \end{split}
\end{align}
and by Lemma \ref{Derivative comutation} and \eqref{bij-xi},
\begin{align}\label {III112}
\begin{split}
III_{112}  =& 2u_n ^{l + 3} \sigma _l \left( G \right)\sum\limits_{j \in B} {u_{njj} \left( {\Delta u - u_{nn} } \right)}  - 2u_n ^{l + 3} \sigma _l \left( G \right)\sum\limits_{i = 1}^{n - 1} {\sum\limits_{j \in B} {u_{ni} u_{ijj} } }  \\
  = & - 2u_n ^{l + 4} \sigma _l \left( G \right)\sigma _1 \left( G \right)\sum\limits_{j \in B} {\left( {u_{jjn}  + \epsilon u_n } \right)}  - 2u_n ^{l + 3} \sigma _l \left( G \right)\sum\limits_{i = 1}^{n - 1} {\sum\limits_{j \in B} {u_{ni} u_{jji} } }  \\
  = & - 4u_n ^{l + 3} \sigma _l \left( G \right)\sigma _1 \left( G \right)\sum\limits_{j \in B} {u_{nj} ^2 }  - 2\left( {n - l-1} \right)\epsilon u_n ^{l + 5} \sigma _l \left( G \right)\sigma _1 \left( G \right) \ .
\end{split}
\end{align}

Therefore, combining \eqref{III11}, \eqref{III111} and \eqref{III112} we get
\begin{align}\label {III11Over}
\begin{split}
III_{11}  =  - 6u_n ^{l + 3} \sigma _l \left( G \right)\sigma _1 \left( G \right)\sum\limits_{j \in B} {u_{nj} ^2 }  - 2\left( {n - l-1} \right)\epsilon u_n ^{l + 5} \sigma _l \left( G \right)\sigma _1 \left( G \right).
\end{split}
\end{align}

For the term $III_{12}$, we have by \eqref{bij-xi} that
\begin{align}\label {III12Over}
\begin{split}
III_{12}  =& 2u_n ^{l + 2} \sigma _l \left( G \right)\sum\limits_{\alpha ,\beta } {\sum\limits_{j \in B} {\left( {a^{\alpha \beta } } \right)_j u_{\alpha \beta j} } }  \\
  = &4u_n ^{l + 3} \sigma _l \left( G \right)\sum\limits_\alpha  {\sum\limits_{j \in B} {u_{nj} u_{\alpha \alpha j} } }  - 4u_n ^{l + 3} \sigma _l \left( G \right)\sum\limits_\alpha  {\sum\limits_{j \in B} {u_{\alpha j} u_{\alpha nj} } }  \\
  =& 4u_n ^{l + 3} \sigma _l \left( G \right)\sum\limits_{i = 1}^{n - 1} {\sum\limits_{j \in B} {u_{nj} u_{iij} } }  \\
  = & - 4u_n ^{l + 4} \sigma _l \left( G \right)\sum\limits_{i \in G,j \in B} {u_{nj} b_{ii,j} }  - 4u_n ^{l + 3} \sigma _l \left( G \right)\sigma _1 \left( G \right)\sum\limits_{j \in B} {u_{nj} ^2 }  \ .
\end{split}
\end{align}

Therefore, by \eqref{III1},  \eqref{III11Over} and \eqref{III12Over} we have
\begin{align}\label {III1Over}
\begin{split}
 III_1  =&  - u_n ^{l + 2} \sigma _l \left( G \right)\sum\limits_{j \in B} {f_{jj} }  - 4u_n ^{l + 4} \sigma _l \left( G \right)\sum\limits_{i \in G,j \in B} {u_{nj} b_{ii,j} }  \\
  &- 10u_n ^{l + 3} \sigma _l \left( G \right)\sigma _1 \left( G \right)\sum\limits_{j \in B} {u_{nj} ^2 }  - 2\left( {n - l-1} \right)\epsilon u_n ^{l + 5} \sigma _l \left( G \right)\sigma _1 \left( G \right) \ .
\end{split}
\end{align}

So, by \eqref{III}, \eqref{III2+III3} and \eqref{III1Over}, it follows that
\begin{align}\label {IIIOver}
\begin{split}
III   \sim & - u_n ^{l + 2} \sigma _l \left( G \right)\sum\limits_{j \in B} {f_{jj} }  - 4u_n ^{l + 4} \sigma _l \left( G \right)\sum\limits_{i \in G,j \in B} {u_{nj} b_{ii,j} }  \\
  &- 10u_n ^{l + 3} \sigma _l \left( G \right)\sigma _1 \left( G \right)\sum\limits_{j \in B} {u_{nj} ^2 }  - 2\left( {n - l-1} \right)\epsilon u_n ^{l + 5} \sigma _l \left( G \right)\sigma _1 \left( G \right) \\
 &+2\left( {n - l-1} \right)\epsilon fu_n ^{l + 2} \sigma _l \left( G \right) \ .
\end{split}
\end{align}

For the term $IV$, we show by Lemma \ref{Derivative comutation} and \eqref{daoshuzhuanyi} that
\begin{align}\label {IVOver}
\begin{split}
 IV =& 2( - 1)^{l + 1} \sum\limits_{\alpha ,\beta ,\xi ,\eta } {a^{\xi \eta } \frac{{\partial \sigma _{l + 2} (D^2 u)}}{{\partial u_{\alpha \beta } }}u_{\alpha \xi \eta } u_\beta  }  \\
  = &2( - 1)^{l + 1} u_n \sum\limits_{\xi ,\eta } {a^{\xi \eta } \frac{{\partial \sigma _{l + 2} (D^2 u)}}{{\partial u_{nn} }}u_{n\xi \eta } }  + 2( - 1)^{l + 1} u_n \sum\limits_{\xi ,\eta } {a^{\xi \eta } \sum\limits_{i = 1}^{n - 1} {\frac{{\partial \sigma _{l + 2} (D^2 u)}}{{\partial u_{in} }}u_{i\xi \eta } } }  \\
  =& 2( - 1)^{l + 1} u_n \sum\limits_{\xi ,\eta } {a^{\xi \eta } \sigma _{l + 1} (u_{pq} )u_{n\xi \eta } }  + 2( - 1)^l u_n \sum\limits_{\xi ,\eta } {a^{\xi \eta } \sum\limits_{i = 1}^{n - 1} {\sigma _l (u_{pq} \left| i \right.)u_{ni} u_{i\xi \eta } } }  \\
  \sim& 2u_n ^{l + 1} \sigma _l (G)\sum\limits_{\xi ,\eta } {a^{\xi \eta } \sum\limits_{j \in B} {u_{nj} u_{j\xi \eta } } }  \\
  =& 2u_n ^{l + 1} \sigma _l (G)\sum\limits_{j \in B} {u_{nj} \sum\limits_{\xi ,\eta } {a^{\xi \eta } } u_{\xi \eta j} }  \\
  = &2u_n ^{l + 1} \sigma _l (G)\sum\limits_{j \in B} {u_{nj} f_j }  + 4u_n ^{l + 3} \sigma _l (G)\sigma _1 (G)\sum\limits_{j \in B} {u_{nj} ^2 }  \ .
 \end{split}
\end{align}

For the term $V$, we have
\begin{align}\label {V}
\begin{split}
 V = &2( - 1)^{l + 1} \sum\limits_{\alpha ,\beta ,\xi ,\eta } {a^{\xi \eta } \frac{{\partial \sigma _{l + 2} (D^2 u)}}{{\partial u_{\alpha \beta } }}u_{\alpha \xi } u_{\beta \eta } }  \\
  = &2( - 1)^{l + 1} \sum\limits_{\xi ,\eta } {a^{\xi \eta } \frac{{\partial \sigma _{l + 2} (D^2 u)}}{{\partial u_{nn} }}u_{n\xi } u_{n\eta } }  + 2( - 1)^{l + 1} \sum\limits_{\xi ,\eta } {a^{\xi \eta } \sum\limits_{i,j = 1}^{n - 1} {\frac{{\partial \sigma _{l + 2} (D^2 u)}}{{\partial u_{ij} }}u_{i\xi } u_{j\eta } } }  \\
  &+ 4( - 1)^{l + 1} \sum\limits_{\xi ,\eta } {a^{\xi \eta } \sum\limits_{j = 1}^{n - 1} {\frac{{\partial \sigma _{l + 2} (D^2 u)}}{{\partial u_{nj} }}u_{n\xi } u_{j\eta } } }  \\
  =&V_1+V_2+V_3\ .
 \end{split}
\end{align}

It is obvious to see that
\begin{align}\label {V1Over}
\begin{split}
V_1\sim0\ ,
\end{split}
\end{align}
and
\begin{align}\label {V3Over}
\begin{split}
 V_3  =& 4( - 1)^{l + 1} \sum\limits_{\xi ,\eta } {a^{\xi \eta } \sum\limits_{j = 1}^{n - 1} {\frac{{\partial \sigma _{l + 2} (D^2 u)}}{{\partial u_{nj} }}u_{n\xi } u_{j\eta } } }  \\
  \sim& 4( - 1)^l \sum\limits_\xi  {a^{\xi \xi } \sum\limits_{j \in B}^{} {\sigma _l (u_{pq} \left| j \right.)u_{jn} u_{n\xi } u_{j\xi } } }  \\
  \sim& 4u_n ^l u_{nn} \sigma _l (G)\sum\limits_{j \in B} {u_{jn} ^2 }  \ .
 \end{split}
\end{align}

For the term $V_2$, we have
\begin{align}\label {V2}
\begin{split}
 V_2  =& 2( - 1)^{l + 1} \sum\limits_{\xi ,\eta } {a^{\xi \eta } \sum\limits_{i,j = 1}^{n - 1} {\frac{{\partial \sigma _{l + 2} (D^2 u)}}{{\partial u_{ij} }}u_{i\xi } u_{j\eta } } }  \\
  =& 2( - 1)^{l + 1} \sum\limits_\xi  {a^{\xi \xi } \sum\limits_{i = 1}^{n - 1} {\frac{{\partial \sigma _{l + 2} (D^2 u)}}{{\partial u_{ii} }}u_{i\xi } ^2 } }  + 2( - 1)^{l + 1} \sum\limits_\xi  {a^{\xi \xi } \sum\limits_{i,j = 1,i \ne j}^{n - 1} {\frac{{\partial \sigma _{l + 2} (D^2 u)}}{{\partial u_{ij} }}u_{i\xi } u_{j\xi } } }  \\
  =& V_{21} +V_{22}  \ .
 \end{split}
\end{align}

For the term $V_{21}$, by \eqref{Formula},
\begin{align}\label {V21}
\begin{split}
V_{21}  =& 2( - 1)^{l + 1} \sum\limits_\xi  {a^{\xi \xi } \sum\limits_{i = 1}^{n - 1} {\frac{{\partial \sigma _{l + 2} (D^2 u)}}{{\partial u_{ii} }}u_{i\xi } ^2 } }  \\
 =&2( - 1)^{l + 1} \sum\limits_\xi  {a^{\xi \xi } \sum\limits_{i = 1}^{n - 1} {\left( {u_{nn} \sigma _l \left( {u_{pq} \left| i \right.} \right) - \sum\limits_{j \ne i} {u_{nj} ^2 \sigma _{l - 1} \left( {u_{pq} \left| {ij} \right.} \right)} } \right)u_{i\xi } ^2 } }  \\
  \sim&  - 2u_n ^l u_{nn} \sigma _l \left( G \right)\sum\limits_{j \in B} {u_{nj} ^2 }  - 2u_n ^{l - 1} \sum\limits_\xi  {a^{\xi \xi } \sum\limits_{i = 1}^{n - 1} {\sum\limits_{j \ne i} {\sigma _{l - 1} \left( {\lambda \left| {ij} \right.} \right)u_{nj} ^2 u_{i\xi } ^2 } } }  \ .
\end{split}
\end{align}

As the calculations before, we divide the term including $i\neq j$ into four parts and then derive
\begin{align}\label {V21Over}
\begin{split}
V_{21} \sim & - 2u_n ^l u_{nn} \sigma _l \left( G \right)\sum\limits_{j \in B} {u_{nj} ^2 }  - 2u_n ^{l - 1} \sum\limits_\xi  {a^{\xi \xi } \sum\limits_{i \in G,j \in B} {\sigma _{l - 1} \left( {G\left| i \right.} \right)u_{nj} ^2 u_{i\xi } ^2 } }  \\
  &- 2u_n ^{l - 1} \sum\limits_{i \in B,j \in G} {\sigma _{l - 1} \left( {G\left| j \right.} \right)u_{nj} ^2 u_{ni} ^2 }  - 2u_n ^{l - 1} \sigma _{l - 1} \left( G \right)\sum\limits_{i,j \in B,i \ne j} {u_{nj} ^2 u_{ni} ^2 }  \\
 \sim & - 2u_n ^l u_{nn} \sigma _l \left( G \right)\sum\limits_{j \in B} {u_{nj} ^2 }  - 4u_n ^{l - 1} \sum\limits_{i \in G,j \in B} {\sigma _{l - 1} \left( {G\left| i \right.} \right)u_{nj} ^2 u_{ni} ^2 }  \\
  &- 2u_n ^{l + 1} \left( {1 + u_n ^2 } \right)\sigma _l \left( G \right)\sigma _1 \left( G \right)\sum\limits_{j \in B} {u_{nj} ^2 }  - 2u_n ^{l - 1} \sigma _{l - 1} \left( G \right)\sum\limits_{i,j \in B,i \ne j} {u_{nj} ^2 u_{ni} ^2 }  \ .
\end{split}
\end{align}

For the term $V_{22}$, it follows that
\begin{align}\label {V22Over}
\begin{split}
 V_{2{\rm{2}}}  =& 2( - 1)^{l + 1} \sum\limits_\xi  {a^{\xi \xi } \sum\limits_{i \ne j}^{} {\frac{{\partial \sigma _{l + 2} (D^2 u)}}{{\partial u_{ij} }}u_{i\xi } u_{j\xi } } }  \\
  = &2( - 1)^{l + 1} \sum\limits_\xi  {a^{\xi \xi } \sum\limits_{i \ne j}^{} {u_{ni} u_{nj} \sigma _{l - 1} \left( {u_{pq} \left| {ij} \right.} \right)u_{i\xi } u_{j\xi } } }  \\
  \sim& 4u_n ^{l - 1} \sum\limits_{i \in G,j \in B}^{} {\sigma _{l - 1} \left( {G\left| i \right.} \right)u_{ni} ^2 u_{nj} ^2 }  + 2u_n ^{l - 1} \sigma _{l - 1} \left( G \right)\sum\limits_{i,j \in B,i \ne j}^{} {u_{ni} ^2 u_{nj} ^2 }  \ .
 \end{split}
\end{align}

By \eqref{V2}, \eqref{V21Over} and \eqref{V22Over},

\begin{align}\label {V2Over}
\begin{split}
V_2   \sim  - 2u_n ^l u_{nn} \sigma _l \left( G \right)\sum\limits_{j \in B} {u_{nj} ^2 }  - 2u_n ^{l + 1} \left( {1 + u_n ^2 } \right)\sigma _l \left( G \right)\sigma _1 \left( G \right)\sum\limits_{j \in B} {u_{nj} ^2 }\ .
\end{split}
\end{align}

Therefore, by \eqref{V}, \eqref{V1Over}, \eqref{V3Over} and \eqref{V3Over} we derive that
\begin{align}\label {VOver}
\begin{split}
V   \sim &2u_n ^l u_{nn} \sigma _l \left( G \right)\sum\limits_{j \in B} {u_{nj} ^2 }  - 2u_n ^{l + 1} \left( {1 + u_n ^2 } \right)\sigma _l \left( G \right)\sigma _1 \left( G \right)\sum\limits_{j \in B} {u_{nj} ^2 }  \\
  =& 2u_n ^l f\sigma _l \left( G \right)\sum\limits_{j \in B} {u_{nj} ^2 }  \ .
\end{split}
\end{align}

Now, combine \eqref{Expression}, \eqref{I+II}, \eqref{IIIOver}, \eqref{IVOver} and \eqref{VOver} we derive finally that
\begin{align}\label {ExpressionOver}
\begin{split}
\sum\limits_{\alpha,\beta}a^{\alpha\beta}\varphi_{\alpha\beta}\sim&- 2u_n ^{l + 3} \sum\limits_\xi  {a^{\xi \xi } \sum\limits_{i \in G,j \in B} {\sigma _{l - 1} \left( {G\left| i \right.} \right)b_{ij,\xi } ^2 } }  - u_n ^{l + 3} \sigma _{l - 1} \left( G \right)\sum\limits_{i,j \in B} {b_{ij,n} ^2 }  \\
  &- u_n ^{l + 3} \left( {1 + u_n ^2 } \right)\sigma _{l - 1} \left( G \right)\sum\limits_{k = 1}^{n - 1} {\sum\limits_{i,j \in B} {b_{ij,k} ^2 } }  + 4u_n ^{l + 2} \sigma _l \left( G \right)\sum\limits_{i \in G,j \in B} {u_{nj} b_{ii,j} }  \\
  &- 2u_n ^{l + 1} \sigma _l \left( G \right)\sigma _1 \left( G \right)\sum\limits_{j \in B} {u_{nj} ^2 }  - 2\left( {n - l-1} \right)\epsilon u_n ^{l + 5} \sigma _l \left( G \right)\sigma _1 \left( G \right) \\
  &- u_n ^{l + 2} \sigma _l \left( G \right)\sum\limits_{j \in B} {f_{jj} }  - 6fu_n ^l \sigma _l \left( G \right)\sum\limits_{j \in B} {u_{nj} ^2 }  + 6u_n ^{l + 1} \sigma _l \left( G \right)\sum\limits_{j \in B} {u_{jn} f_j }  \\
 &+2\left( {n - l-1} \right)\epsilon fu_n ^{l + 2} \sigma _l \left( G \right) \ .
\end{split}
\end{align}

Since $a^{\xi\xi}\ge 1$ for any $\xi=1,2,\cdots,n$, $b_{ij,i}=b_{ii,j}$ for $i\in G,\, j\in B$ and $\epsilon\ge 0$ , we then have
\begin{align}\label {ExpressionInequality}
\begin{split}
\sum\limits_{\alpha,\beta}a^{\alpha\beta}\varphi_{\alpha\beta} \preceq &- 2u_n ^{l + 3}  \sum\limits_{i \in G,j \in B} {\sigma _{l - 1} \left( {G\left| i \right.} \right)b_{ii,j } ^2 } + 4u_n ^{l + 2} \sigma _l \left( G \right)\sum\limits_{i \in G,j \in B} {u_{nj} b_{ii,j} }  \\
   &- 2u_n ^{l + 1} \sigma _l \left( G \right)\sigma _1 \left( G \right)\sum\limits_{j \in B} {u_{nj} ^2 }   \\
  &- u_n ^{l + 2} \sigma _l \left( G \right)\sum\limits_{j \in B} {f_{jj} }  - 6fu_n ^l \sigma _l \left( G \right)\sum\limits_{j \in B} {u_{nj} ^2 }  + 6u_n ^{l + 1} \sigma _l \left( G \right)\sum\limits_{j \in B} {u_{jn} f_j }  \\
 &+2\left( {n - l-1} \right)\epsilon fu_n ^{l + 2} \sigma _l \left( G \right) \ .
\end{split}
\end{align}

After a observation, we can simplify the above expression into the following form
\begin{align}\label {Ercihanshu}
\begin{split}
\sum\limits_{\alpha,\beta}a^{\alpha\beta}\varphi_{\alpha\beta} \preceq &- 2u_n ^{l +1}  \sum\limits_{i \in G,j \in B}\left(u_n\sqrt{\sigma_{l-1}(G|i)}b_{ii,j }-
\sqrt{\sigma_l(G)\lambda_i}u_{nj} \right)^2 \\
  &- u_n ^{l } \sigma _l \left( G \right)\sum\limits_{j \in B}\left[ {6f {u_{nj} ^2 }  -6u_n f_j u_{n j}+ u_n^2(f_{jj}-2\epsilon f)    } \right]\ .
\end{split}
\end{align}

Now, we come to compute the term involving $ f$. Note that $f(x,\,\nabla u)=H(x)(1+|\nabla u|^2)^{\frac 3 2}$, then we get for $j\in B$,
\begin{align}\label {f01}
\begin{split}
f(x,\,\nabla u)&=H(1+u_n^2)^{\frac 3 2}\,,\\
\left[f(x,\,\nabla u)\right]_j&=H_j(1+u_n^2)^{\frac 3 2}+3Hu_n(1+u_n^2)^{\frac 1 2}u_{nj}\,,
\end{split}
\end{align}
and
\begin{align}\label {f2}
\begin{split}
\left[f(x,\,\nabla u)\right]_{jj}=&H_{jj}(1+u_n^2)^{\frac 3 2}+6H_ju_n(1+u_n^2)^{\frac 1 2}u_{nj}+3Hu_n^2(1+u_n^2)^{-\frac 1 2}u_{nj}^2\\
&+3H(1+u_n^2)^{\frac 1 2}u_{nj}^2+3H(1+u_n^2)^{\frac 1 2}u_n u_{njj}\,.
\end{split}
\end{align}

By Lemma \ref{Derivative comutation}, \eqref{bij-xi} and the fact that $\epsilon\ge 0$, we have for $j\in B$,
\begin{align}\label {f21}
\begin{split}
\left[f(x,\,\nabla u)\right]_{jj}=&H_{jj}(1+u_n^2)^{\frac 3 2}+6H_ju_n(1+u_n^2)^{\frac 1 2}u_{nj}+3Hu_n^2(1+u_n^2)^{-\frac 1 2}u_{nj}^2\\
&+9H(1+u_n^2)^{\frac 1 2}u_{nj}^2+3Hu_n^2(1+u_n^2)^{\frac 1 2}\epsilon-3Hu_n^2(1+u_n^2)^{\frac 1 2}b_{jj,n}\\
\ge &H_{jj}(1+u_n^2)^{\frac 3 2}+6H_ju_n(1+u_n^2)^{\frac 1 2}u_{nj}+3Hu_n^2(1+u_n^2)^{-\frac 1 2}u_{nj}^2\\
&+9H(1+u_n^2)^{\frac 1 2}u_{nj}^2-3Hu_n^2(1+u_n^2)^{\frac 1 2}b_{jj,n}\,.
\end{split}
\end{align}

Therefore, by \eqref{1-ordercondition} we derive
\begin{align}\label {f21}
\begin{split}
\sum_{j\in B}\left[6f {u_{nj} ^2 }  -6u_n f_j u_{n j}+ u_n^2(f_{jj}-2\epsilon f)\right] \ge \sum_{j\in B} \left(Au_{nj}^2+Bu_{nj}+C \right)\,,
\end{split}
\end{align}
where
\begin{align}\label {f21ABC}
\begin{split}
A=&3H(2+u_n^2)(1+u_n^2)^{-\frac 1 2}\ge 6H(1+u_n^2)^{-\frac 1 2}\,,\\
B=&-6u_nH_j(1+u_n^2)^{\frac 1 2}\,,\\
C=&u_n^2(H_{jj}-2\epsilon H)(1+u_n^2)^{\frac 3 2}\,.
\end{split}
\end{align}

Under the structure condition \ref{structurecondition}, we easily get that for $j\in B$,
\begin{align}\label {f22}
\begin{split}
 Au_{nj}^2+Bu_{nj}+C\ge 0\,.
\end{split}
\end{align}

Thus,
\begin{align}\label {final}
\begin{split}
\sum\limits_{\alpha,\beta}a^{\alpha\beta}\varphi_{\alpha\beta} \preceq 0.
\end{split}
\end{align}
And this finishes the proof of the constant rank theorem.\hfill $\sharp$

\section{Strict Convexity of Level sets}
\setcounter{equation}{0} \setcounter{theorem}{0}\noindent

In this section, we follow the idea of \cite{Korevaar} to prove the remained parts of Theorem \ref{WZ} by the continuity method and then to prove  Corollary \ref{WZCor}.

Remark that the equation we concentrate is the following
\begin{align}\label {Equation-u-tau}
\begin{split}
\left\{ {\begin{array}{lll}
\mathrm{div}(\frac{\nabla u^\tau}{\sqrt{1+|\nabla u^\tau|^2}}) =0 & \mathrm{\mbox {in}} & \ \Omega=\Omega_0\setminus\bar{\Omega}_1,\\
u^\tau=0& \mathrm{\mbox {on}}& \  \partial \Omega_0,\\
u^\tau=\tau& \mathrm{\mbox {on}}&  \ \partial \Omega_1
\end{array}} \right.
\end{split}
\end{align}
for $0<\tau\le 1.$

Firstly, when the height of the minimal graph, $\tau$, is small enough, the level sets must be regular and strictly convex. In fact, if $\tau $ is small, by Theorem \ref{WZtau} we know that the norm of the gradient will be sufficiently small which would force the minimal graph to be looked as something like the graph of a harmonic function since the similarity of the equations they satisfy. Precisely, we have
\begin{lemma} \label{smallheight}
Let $(M^n,g)$ be a space form with constant sectional curvature $\epsilon\ge0$ , and $\Omega_0$ and $\Omega_1$ be bounded smooth strictly convex domains in $M^n$, $n\ge 2$ and $\bar{\Omega}_1\subseteq \Omega_0$. Then there exists $\delta_0>0$ such that for any $0<\tau\le \delta_0$, the solution to the minimal graph equation \eqref{Equation-u-tau} satisfies that $\nabla u^\tau \neq 0$ and all level sets of $u^\tau$ are strictly convex with respect to $\nabla u^\tau$.
\end{lemma}

{\bf Proof:\ \ }
Let $\omega^\tau$ be the harmonic function defined in \eqref{omegatau}. Denote by $\rho=u^\tau-\omega^\tau$ and it is easy to see that $\rho$ satisfies the following possion type equation
 \begin{align}\label {rho}
\begin{split}
\left\{ {\begin{array}{lll}
\triangle \rho=\sum\limits_{\xi,\eta}u^\tau_{\xi\eta}u^\tau_\xi u^\tau_\eta-|\nabla u^\tau|^2\triangle u^\tau & \mathrm{\mbox {in}} & \ \Omega=\Omega_0\setminus\bar{\Omega}_1,\\
\rho=0& \mathrm{\mbox {on}}& \  \partial \Omega_0,\\
\rho=0& \mathrm{\mbox {on}}&  \ \partial \Omega_1.
\end{array}} \right.
\end{split}
\end{align}

 Denoted by $f=\sum\limits_{\xi,\eta}u^\tau_{\xi\eta}u^\tau_\xi u^\tau_\eta-|\nabla u^\tau|^2\triangle u^\tau$, according to Theorem \ref{WZtau} we know that
 \begin{align}\label {fupperbound}
\begin{split}
|f|_{C^{\alpha}(\bar{\Omega})}\le C \tau^2.
\end{split}
\end{align}

 The Schauder theory then gives that
\begin{align}\label {rhoestimate}
\begin{split}
||u^\tau-\omega^\tau||_{C^{2,\alpha}(\bar{\Omega})}\le & C \left( ||\rho||_{L^\infty(\Omega)}+|f|_{C^\alpha(\bar{\Omega})}    \right)\\
\le & C \left( ||f||_{L^\infty(\bar{\Omega})}+|f|_{C^\alpha(\bar{\Omega})}    \right)
 \le C_5\tau^2.
\end{split}
\end{align}

Especially, we get for any $x\in \Omega$,
\begin{align}\label {gradient-estimate}
\begin{split}
|\nabla u^\tau-\nabla \omega^\tau|(x)\le  C_5\tau^2.
\end{split}
\end{align}
It follows by \eqref{lip0} that for any $x\in \Omega$,
\begin{align}
\begin{split}
|\nabla u^\tau|(x)\ge | \nabla \omega^\tau |(x) -C_5\tau^2\ge C_0\tau-C_5\tau^2,
\end{split}
\end{align}
thus there exists some $\delta_0>0$ such that for any $\tau\in (0,\, \delta_0]$ and for any $x\in \Omega$,
\begin{align}\label {u-tau-estimate}
\begin{split}
|\nabla u^\tau|(x)> 0\,,
\end{split}
\end{align}
which guarantees the regularity of the level sets in this situation.

Next, we come to consider the strict convexity of the level sets.

 Notations as before, according to \cite{Mzyb} we know that the level sets of $\omega^1$ is strictly convex, then
\begin{align}\label {eigenvalue}
\begin{split}
\lambda_{\mathrm{min}}(\omega^1)\ge C_6
\end{split}
\end{align}
for some positive constant $C_6$, where $\lambda_{\mathrm{min}}(\psi)$ is denoted to be the minimal eigenvalue of the second fundamental forms of the level sets of $\psi$ with respect to the gradient direction. Remark that $\lambda_{\mathrm{min}}(\psi)$ is an expression via the derivatives up to $2-$order of $\psi$.

Since the relation of the level sets between $\omega^\tau$ and $\omega^1$ is that for $s\in [0,\,\tau ]$
$$(\omega^\tau)^{-1} (s)=(\omega^1)^{-1}(\frac{s}{\tau}),$$
we deduce that for any $0<\tau\le 1$
\begin{align}\label {eigenvalue-tau}
\begin{split}
\lambda_{\mathrm{min}}(\omega^\tau)\ge C_6.
\end{split}
\end{align}

By \eqref{rhoestimate}, $u^\tau$ and $\omega^\tau$ is $C^{2,\alpha}$ close enough once $\tau$ is small enough, so by adjusting suitably $\delta_0$ above, we can get that for $\tau\in (0,\, \delta_0]$,
$$\lambda_{\mathrm{min}}(u^\tau)\ge C_7>0.$$
This gives the strict convexity of the level sets of minimal graph with a small height.\hfill $\sharp$

{\bf Proof of Theorem \ref{WZ}:}\ \ Since the existence of minimal graph defined on a convex ring is already settled in Theorem \ref{WZtau}, we now begin to conclude the remained results.

Set
\begin{align*}
T=\{\tau\in (0,\,1]\,|\, \forall s\, \in (0,\tau], u^s \mbox{ has nonzero gradient and stirctly convex level sets }   \}.
\end{align*}

By Lemma \ref{smallheight}, we know that $(0,\ \delta_0]\subset T.$

For the openness, we assume firstly that $\tau_0\in T$, that is to say, for all $x\in \Omega$ we have $|\nabla u^{\tau_0}|(x)>0$ and the level sets of $u^{\tau_0}$ has positive second fundamental form with respect to its gradient. Namely, there is a $C_8>0$ such that
 \begin{align}\label {open-1}
\begin{split}
|\nabla u^{\tau_0}|\ge C_8,\ \lambda_{\mathrm{min}}(u^{\tau_0})\ge C_8.
\end{split}
\end{align}
 Then Theorem \ref{WZtau} and almost the same procedure as Lemma \ref{smallheight} produce that there exists some $\delta>0$ such that for all $\tau \in [\tau_0-\delta,\ \tau_0+\delta ]$, $u^\tau$ has nonsingular gradient and strictly convex level sets, this concludes the openness.

For the closedness, we assume without loss of generality that
$$\tau_i\in \ T,\ \tau_\infty\in(0,\ 1],\ \delta_0\le\tau_i<\tau_\infty,\ \ \mbox{and}\ \lim_{i\rightarrow \infty}\tau_i=\tau_\infty.$$
Obviously, $u^{\tau_\infty}$ is the unique solution to \eqref{Equation-u-tau} at the time $\tau_\infty $. In the sequel, we will prove that $\tau_\infty \in T,$ that is to say, $u^{\tau_\infty}$ has nonzero gradient and strictly convex level sets.

For each $\tau_i \in T$, according to the assumption we have $|\nabla u^{\tau_i}|>0$ and the level sets of $u^{\tau_i}$ are strictly convex with respect to $\nabla u^{\tau_i}$. Discussion as before, taking the frame such that $e_n=\frac{\nabla u^{\tau_i}}{|\nabla u^{\tau_i}|}$, we can derive by \eqref{BBB} that $ u^{\tau_i}_{ll}<0$ for $l=1,\,2,\,\cdots,n-1$. Thus we deduce from the equation of minimal graph that
\begin{align}\label {unn}
\begin{split}
\frac{{\sum\limits_{\alpha ,\beta  = 1}^n {u^{\tau_i} _{\alpha \beta } u^{\tau_i} _\alpha  u^{\tau_i} _\beta  } }}{{|\nabla u^{\tau_i} |^2 }} = u^{\tau_i} _{nn}  =  -(1+u_n^2) \sum\limits_{l = 1}^{n - 1} {u^{\tau_i} _{ll} }   > 0.
\end{split}
\end{align}
Therefore, the norm of the gradient of $u^{\tau_i}$ will increase along the gradient direction, in fact,
\begin{align}\label {Decrease}
\begin{split}
\sum\limits_{\alpha=1}^{n}  u^{\tau_i}_\alpha|\nabla u^{\tau_i}|^2_\alpha =2\sum\limits_{\alpha ,\beta  = 1}^n {u^{\tau_i} _{\alpha \beta } u^{\tau_i} _\alpha  u^{\tau_i} _\beta  } >0.
\end{split}
\end{align}
Thus the minimum of $|\nabla u^{\tau_i}|$ has to be attained on the exterior boundary $\partial \Omega_0.$

On the other hand, it is obvious from comparison principle that for any $i$,
$$u^{\delta_0}\le  u^{\tau_i}\ \ \mbox{and}\ \ u^{\delta_0}=  u^{\tau_i}=0 \ \ \mbox{on}\ \  \partial \Omega_0,$$
therefore according to Hopf strong maximum principle, a positive constant $C_{10}$ depending upon $u^{\delta_0}$ can be taken to assure that for all $x\in \partial \Omega_0$,
 \begin{align}\label {C10}
\begin{split}
|\nabla u^{\tau_i}|(x)>|\nabla u^{\delta_0}|(x)\ge C_{10}>0.
\end{split}
\end{align}
Thus,
\begin{align}\label {C101}
\begin{split}
\forall x\in  \Omega,\ \ |\nabla u^{\tau_i}|(x)\ge C_{10}>0.
\end{split}
\end{align}

Once again taking advantage of Theorem \ref{WZtau} and the fact that $\tau_i\rightarrow \tau_\infty$, we then have
\begin{align}\label {GradientReguler}
\begin{split}
|\nabla u^{\tau_\infty}|>0.
\end{split}
\end{align}

Finally, for the strict convexity of the level sets of each $u^{\tau_i}$ and Theorem \ref{WZtau}, the limit $u^{\tau_\infty}$ could only has convex level sets. Then the constant rank Theorem \ref{WZ3} and the strict convexity of $\left(u^{\tau_\infty}\right)^{-1}(0)=\partial \Omega_0$ ensure that the level sets of $u^{\tau_\infty}$ must be all strictly convex. This is to say that $\tau_\infty\in T$ then $T$ is relatively closed in $(0,\, 1].$

In summary, we have by the connection of $(0,\,1]$ that $T=(0,\, 1]$, and this finally concludes the whole proof of Theorem \ref{WZ}.\hfill $\sharp$

{\bf Proof of Corollary \ref{WZCor}:}\ \ Based on Theorem \ref{WZ}, we know that the level sets of $u$ are regular and strictly convex. Then the corollary is immediately derived by the same discussion as \eqref{unn} and \eqref{Decrease}.\hfill $\sharp$

 {\bf Acknowledgments:} The research was supported by NSFC (No.11471188). Both authors would like to owe thanks to Prof. X. Ma for his constant encouragement and they also would like to owe thanks to Prof. Chuanqiang Chen for his helpful discussion and suggestions to revise the paper.

\end{document}